\documentstyle[11pt,leqno,amscd,amssymb,latexsym,amsbsy,xypic,mathrsfs,pstricks,color,verbatim]{amsart}
\include{micro}
\def\da{\delta}
\def\gz{\gamma}
\def\ggz{\Gamma}
\def\az{\alpha}

\def\bz{\beta}
\def\ra{\rightarrow}
\def\lra{\longrightarrow}
\def\ma{{\mathcal{A}}}
\def\mc{{\mathcal{C}}}
\def\ms{\mathcal{S}}
\def\mod{\mbox{\rm mod\,}}
\def\Hom{\mbox{\rm Hom}}
\def\End{\mbox{\rm End}}
\def\Ext{\mbox{\rm Ext}}
\def\rad{\mbox{\rm rad\,}}
\def\add{\mbox{\rm add\,}}

\def\vz{\varepsilon}

\def\t{\tau}

\newtheorem{thm}{Theorem}[section]
\newtheorem{lem}[thm]{Lemma}
\newtheorem{cor}[thm]{Corollary}
\newtheorem{prop}[thm]{Proposition}
\newtheorem{rem}[thm]{Remark}

\theoremstyle{definition}
\newtheorem{example}[thm]{Example}

\newtheorem{defn}[thm]{Definition}

\begin{document}

\title{On the Cluster Category of a Marked Surface}

\author{Thomas Br\"ustle and Jie Zhang}
\address{D\'{e}partement de Math\'{e}matiques, Universit\'{e} de Sherbrooke, Sherbrooke, Canada, J1K 2R1 and Department of Mathematics, Bishop¡¯s University, Sherbrooke, Canada, J1M 1Z7} \email{thomas.brustle@@usherbrooke.ca and tbruestl@@ubishops.ca}
\address{D\'{e}partement de Math\'{e}matiques, Universit\'{e} de Sherbrooke, Sherbrooke, Canada, J1K 2R1} \email{Jie.Zhang@@usherbrooke.ca}
\thanks{Thomas Br\"ustle is supported by  NSERC.
Jie Zhang is partially supported by China Scholarship Council  and Bishop's University.
Both authors stayed at INRIA Sophia Antipolis when writing this paper.}
\date{\today}

\thanks{}
\begin{abstract}
We study in this paper the cluster category $\mc_{(S,M)}$ of a marked surface $(S,M)$ without punctures. We explicitly describe the objects in $\mc_{(S,M)}$ as direct sums of  homotopy classes of curves in $(S,M)$ and one-parameter families related to closed curves in $(S,M)$.
Moreover, we describe the Auslander-Reiten structure of the category $\mc_{(S,M)}$ in geometric terms and show that the objects without self-extensions in $\mc_{(S,M)}$  correspond to curves in $(S,M)$ without self-intersections. As a consequence, we establish that every rigid indecomposable object is reachable from an initial triangulation.
\end{abstract}
\maketitle

\section{Introduction}

We study in this paper the cluster category of a marked surface:
Consider a compact connected oriented 2-dimensional bordered  Riemann surface $S$ and a finite set of marked points $M$ lying on the boundary of $S$ with at least one marked point on each boundary.
In \cite{FST} a cluster algebra $\ma(S,M)$ is associated to the marked surface $(S,M)$. The initial seed of  $\ma(S,M)$  corresponds to a triangulation $\Gamma$ of $(S,M)$, and the mutation of a cluster variable corresponds to the flip of an arc in the triangulation (see \cite{FZ} for the definition of a cluster algebra). The cluster algebra, defined by iterated mutations, thus is independent of the chosen triangulation $\Gamma$ of $(S,M)$.

The cluster category $\mc_{(S,M)}$ providing a categorification of the algebra  $\ma(S,M)$  has been defined in \cite{A}.
In fact, in \cite{LF,ABCP} a quiver with potential $(Q_{\ggz},W_{\ggz})$ has been defined for each triangulation $\Gamma$ of $(S,M)$.
Since the Jacobian algebra $J(Q_\ggz,W_\ggz)$ is finite-dimensional, one can use \cite{A} to define the cluster category $\mc_\ggz$ associated to $\ggz$.
By \cite{KY,FST,LF}, this category $\mc_\ggz$ does not depend on the triangulation $\ggz$ of $(S,M)$ and is just denoted by $\mc_{(S,M)}$ (more precisely,  $\mc_\ggz$ is triangle equivalent to $\mc_{\ggz'}$ if $\ggz'$ is another triangulation of $(S,M)$).

Cluster categories associated to quivers are object of intense investigation (see for example  \cite{GLS1,CCS,BMRRT,BMR,CC,GLS2,ABS}) and are understood explicitly in terms of representations of the quiver. The category $\mc_{(S,M)}$, defined by a quiver with potential, is however quite difficult to compute and only a few cases of surfaces allow direct calculations. The aim of this paper is to provide an explicit description of the objects and the irreducible morphisms in the category $\mc_{(S,M)}$ in terms of the surface $(S,M)$, independent of the choice of a triangulation. In fact, the category $\mc_{(S,M)}$ is closely related to module categories of string algebras, and a well-know classification theorem \cite{BR} describes all modules over a string algebra by listing the indecomposable modules in a list of string modules and a collection of band modules. Similarly, the indecomposable objects of $\mc_{(S,M)}$ can be described as belonging to two classes which we call string objects and band objects.
The cluster categories $\mc_{(S,M)}$  are $k-$categories where $k$ denotes a  fixed algebraically closed field. Moreover they are Krull-Schmidt, thus it is sufficient to classify indecomposable objects up to isomorphism.

\begin{thm}\label{indec}A parametrization of the isoclasses of indecomposable objects in $\mc_{(S,M)}$ is given by ``string objects" and ``band objects", where
\begin{itemize}
\item[(1)] The string objects are indexed by the homotopy classes of non-contractible curves in $(S,M)$ which are not homotopic to a boundary segment of $(S,M)$, subject to the equivalence relation $\gamma \sim \gamma^{-1}$.
\item[(2)] The band objects are indexed by $k^* \times {\Pi^*_1 (S,M)}/ \sim$, where $k^*=k\setminus\{0\}$ and $\Pi^*_1 (S,M)/\sim$  is given by the non-zero elements of the fundamental group of $(S,M)$ subject to the equivalence relation generated by $a \sim a^{-1}$ and cyclic permutation.
\end{itemize}
\end{thm}

The boundary of the surface $S$ consists of a collection of disjoint circles, each one inheriting an orientation from the orientation of $S$. Given any curve $\gamma$  in $(S,M)$, we denote by $_s\gamma$ the curve obtained from $\gamma$ by moving its starting point clockwise to the next marked point on the boundary. Likewise, moving the ending point of $\gz$ clockwise along the boundary to the next marked point, one obtains a curve $\gamma_e$ from $\gamma$. These moves establish the irreducible morphisms between string objects, as explained in the following theorem. Note that in the special case where $S$ is a disc, this description of irreducible morphisms has been given in \cite{CCS}, establishing one of the first categorifications of a cluster algebra (of type $A$).

\begin{thm}\label{irred}Let $\gz$ be a non-contractible curve in $(S,M)$ which is not homotopic to a boundary segment of $(S,M)$. Then there is an AR-triangle in $\mc_{(S,M)}$ as follows, where $_s\gz$ or $\gz_e$ are zero objects in $\mc_{(S,M)}$ if they are boundary segments:
$$\gz \lra  {_s\gamma} \oplus \gz_e \lra  {_s\gz_e} \lra \gz[1].$$
Moreover, all AR-triangles between string objects in $\mc_{(S,M)}$ are of this form.
\end{thm}

We obtain also the AR-translation in the cluster category $\mc_{(S,M)}$:

\begin{prop}The AR-translation in $\mc_{(S,M)}$ is given by simultaneous counterclockwise rotation of starting and ending point of a curve to the next marked points at the boundary.
\end{prop}

Moreover, from Theorem \ref{irred} one can deduce the shape of the AR-compon\-ents in the category $\mc_{(S,M)}$. In fact, each boundary component of $S$ with $t$ marked points gives rise to a tube of rank $t$ in $\mc_{(S,M)}$. All other components formed by string objects are composed of meshes with exactly two middle terms.

We also study the effect of a change of the triangulation. It is well-known that any two triangulations of $(S,M)$ can be transformed into each other by a sequence of flips $f_i(\ggz)$ which locally change one arc in $\ggz$.
On the other hand, each triangulation yields a cluster-tilting object $T$ in $\mc_{(S,M)}$, and there one can apply mutations $\mu_i(T)$ locally changing one summand. We show that these two operations are compatible:

\begin{thm}Each triangulation $\ggz$ of $(S,M)$ yields a cluster-tilting object $T_\Gamma$ in $\mc_{(S,M)}$, and
$$\mu_i(T_\ggz)={T_{f_i(\ggz)}}.$$
\end{thm}

We then study the relation between extensions in the category
$\mc_{(S,M)}$ and intersections of the corresponding curves. Given any two
curves in $(S,M)$ with intersections, we explicitly construct one or two
new curves, sometimes resolving the intersection, and sometimes increasing
the winding number. These new curves serve as middle term of certain non-split short exact
sequences which allow to prove the following theorem and corollary:

\begin{thm} Curves in $(S,M)$ without self-intersections  correspond to the objects without self-extensions in $\mc_{(S,M)}$.
\end{thm}
\begin{cor}There is a bijection between triangulations of $(S,M)$ and cluster-tilting objects of $\mc_{(S,M)}$. In particular, each indecomposable object without self-extension is reachable from the cluster-tilting object $T_\ggz$ (the initial cluster-tilting object).
\end{cor}

This paper is organized as follows: we first describe in section 2 the indecomposable objects in $\mc_{(S,M)}$, using the description of modules over a string algebra. In section 3 we describe the irreducible morphisms in $\mc_{(S,M)}$ and study the shape of its AR-components. Section 4 is devoted to compare the effect of a flip of an arc in the triangulation and a mutation of the corresponding cluster-tilting object. We also compare it to results obtained for decorated representations of quivers with potential.
Finally, in section 5 we prove that curves  without self-intersections  correspond to objects without self-extensions in $\mc_{(S,M)}$ and we establish the bijection between triangulations and cluster-tilting objects.
\bigskip

Jie Zhang wishes to thank Claire Amiot and Yann Palu for interesting discussions on this problem.

\section{Indecomposable objects in $\mc_{(S,M)}$}

In this section, we choose a triangulation $\ggz$ of $(S,M)$, and by studying the corresponding Jabobian algebra, we give a geometric characterization of the indecomposable objects in $\mc_{(S,M)}$.
By a curve in $(S,M)$, we mean a continuous function $\gz : [0,1] \ra S$ with
$\gz(0),\gz(1)\in M$.
A closed curve is one with $\gamma(0)=\gamma(1)$, and a simple curve is one where $\gamma$ is injective, except possibly at the endpoints.
For a curve $\gz$, we denote by $\gz^{-1}$ the inverse curve $\gz^{-1}: [0,1] \ra S, \; t \mapsto \gz(1-t)$.
We always consider curves up to homotopy, and for any collection of curves we implicitly assume that their mutual intersections are minimal possible in their respective homotopy classes.
We recall from \cite{FST} the definition of a triangulation:
\begin{defn}An arc $\da$ in $(S,M)$ is a simple non-contractible curve in $(S,M)$.
The boundary of $S$ is a disjoint union of cycles, which are subdivided by the points in $M$ into boundary segments. We call an arc $\delta$ a \emph{boundary arc} if it is homotopic
to such a boundary segment.
Otherwise, $\delta$ is said to be an \emph{internal arc}. A \emph{triangulation} of $(S,M)$ is a maximal
collection $\ggz$ of arcs that do not intersect except at their endpoints. We call a triangle $\vartriangle$ in $\ggz$ an \emph{internal triangle} if all edges of $\vartriangle$ are internal arcs.
\end{defn}

\begin{prop}[\cite{FST}]In each triangulation of $(S,M)$, the number
of internal arcs is
$$n = 6g + 3b + c -6$$
where $g$ is the genus of $S$, $b$ is the number of boundary components, and
$c = |M|$ is the number of marked points.
\end{prop}
\bigskip

\subsection{The quiver with potential (Q$_\ggz$,W$_\ggz$)}\label{quiver}

We recall from \cite{ABCP,LF} that each triangulation $\ggz$ of $(S,M)$ yields
 a quiver $Q_\ggz$ with potential $W_\ggz$:

\begin{itemize}
\item[(1)] $Q_{\ggz}=(Q_0,Q_1)$ where the set of vertices $Q_0$ is given by the internal arcs of $\ggz$, and the set of arrows $Q_1$ is defined as follows: whenever there is a triangle $\vartriangle$ in $\ggz$ containing two
internal arcs $a$ and $b$, then there is an arrow $\rho:a \ra b$ in $Q_1$ if $a$ is a predecessor
of $b$ with respect to clockwise orientation at the joint vertex of $a$ and $b$ in $\vartriangle$.

\item[(2)]Every internal triangle $\vartriangle$ in $\ggz$ gives rise to an oriented cycle $\alpha_\vartriangle \beta_\vartriangle\gz_\vartriangle$ in
$Q$, unique up to cyclic permutation of the factors $\alpha_\vartriangle, \beta_\vartriangle, \gz_\vartriangle$. We define
$$W_\ggz =\displaystyle\sum_\vartriangle\alpha_\vartriangle \beta_\vartriangle\gz_\vartriangle $$
where the sum runs over all internal triangles $\vartriangle$ of $\ggz$.
\end{itemize}
Unless we compare two different triangulations, we omit the subscript and denote the quiver with potential defined by $\ggz$ just with $(Q,W)$. We refer to \cite{DWZ} for more details on quivers with potentials, such as the definition of the Jacobian algebra $J(Q,W)$ associated to $(Q,W)$.
From \cite{ABCP} we know that  $J(Q,W)$ is a finite-dimensional string algebra provided $(Q,W)$ is defined by a triangulation of a marked surface as above.

Moreover, we denote by $\Lambda_{(Q,W)}$ the corresponding Ginzburg dg-algebra  (see \cite{G} for more details), and denote by $D(\Lambda_{(Q,W)})$ its derived category (see \cite{K}).
The
\emph{cluster category} $\mc_\ggz={\mc}_{(Q,W)}$ associated to $(Q,W)$ is defined in \cite{A} as the quotient of triangulated categories
$Per\Lambda_{(Q,W)}/D^b(\Lambda_{(Q,W)})$ where $Per\Lambda_{(Q,W)}$ is the thick subcategory of $D(\Lambda_{(Q,W)})$ generated by $\Lambda_{(Q,W)}$ and $D^b(\Lambda_{(Q,W)})$ is the
full subcategory of $D(\Lambda_{(Q,W)})$ of the dg-modules whose homology is of finite total dimension.
\medskip

\begin{thm}[\cite{A},\cite{KZ}]\label{AKZ} If $(Q,W)$ is a quiver with potential whose Jacobian algebra $J(Q,W)$ is finite-dimensional, then
\begin{itemize}
\item[(1)] ${\mc}_\ggz$ is 2-Calabi-Yau, Hom-finite and the image $T_\ggz$ of the free module $\Lambda_{(Q,W)}$ in
the quotient $Per\Lambda_{(Q,W)}/D^b(\Lambda_{(Q,W)})$ is a cluster-tilting object.
\item[(2)] ${\mc}_\ggz/{\mbox{T}_\ggz}$ is equivalent to $\mod J(Q,W)$, the category of finite dimensional modules over $J(Q,W)$. Moreover, the projection functor ${\mc}_\ggz\ra \mod J(Q,W)$ is given by $\Ext^1_{\mc_\ggz}(T_\ggz,-).$
\end{itemize}

\end{thm}

As explained in the introduction, the category ${\mc}_\ggz$ is (up to triangle equivalence) independent of the choice of a triangulation of $(S,M)$, and is therefore denoted by $\mc_{(S,M)}$. Moreover, it is a Krull-Schmidt category, so it is sufficient to describe all indecomposable objects up to isomorphism in order to describe its objects.  From Theorem \ref{AKZ}(2), the  indecomposable objects in $\mc_{(S,M)}$ are either indecomposable modules over the Jacobian algebra $J(Q,W)$ or one of the  $|Q_0|=6g + 3b + |M| -6$ summands of $T_\ggz$. In order to give a description of all indecomposable objects in $\mc_{(S.M)}$, we first study the Jacobian algebra $J(Q,W)$.

\subsection{The string algebra J(Q,W)}
In this subsection, we recall some basic definitions related to string algebras and prove Theorem \ref{indec}.
Recall from \cite{BR} that a finite-dimensional algebra $A$ is a \emph{string algebra} if there is a quiver $Q$ and an admissible ideal $I$ such that  $A=kQ/I$ and the following conditions hold:
\begin{itemize}
\item[(S1)]At each vertex of $Q$ start at most two arrows and stop at most two arrows.
\item[(S2)]For each arrow $\alpha$ there is at most one arrow $\beta$ and at most one arrow
$\da$ such that $\alpha\beta \not\in I $ and $\da\alpha\not\in I.$
\end{itemize}
\medskip

Given an arrow $\bz \in Q_1$, let $s(\bz)$ be its starting point and $e(\bz)$ its ending point.
We denote $\bz^{-1}$ the formal inverse of $\bz$ with $s(\bz^{-1})=e(\bz)$ and $e(\bz^{-1})=s(\bz)$.
A word $w=\az_n\az_{n-1}\cdots\az_1$ of arrows and their formal inverses is called a \emph{string} if $\az_{i+1}\neq \az_i^{-1}, e(\az_i)=s(\az_{i+1})$ for all $1\leq i\leq n-1$, and no subword nor its inverse belongs to $I$. Thus a string $w$ can be viewed as a walk in the quiver $Q$ avoiding the zero relations defining the ideal $I$:
$$w: x_1 \frac{~~~\az_1~}{}x_2 \frac{~~~\az_2~}{}\cdots x_{n-1} \frac{\az_{n-1}}{}x_{n}\frac{~~~\az_n~}{}x_{n+1}$$
where $x_i$ are vertices of $Q$ and $\alpha_i$ are arrows in either direction.
We denote by $s(w)=s(\az_1)$ and $e(w)=e(\az_n)$ the starting point and the ending point of $w$, respectively.
 For technical reasons, we also consider the empty string which we also call zero string.
 A string $w$ is called \emph{cyclic} if the first vertex $x_1$ and the last vertex $x_{n+1}$ coincide.
 A \emph{band} $b=\az_n\az_{n-1}\cdots\az_2\az_1$ is defined to be a
cyclic string $b$ such that each power $b^m$ is a string, but $b$ itself is not a proper
power of any string. Thus $b$ can be viewed as a cyclic walk:

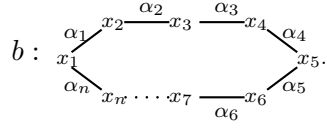
\begin{center}
\begin{pspicture}(-3.4,-.7)(12,0.6)
\put(-.3,0){$b:$}
\put(.3,0){$_{x_1}$} \put(.9,-.5){$_{x_n}$}
\put(1.8,-.5){$_{x_7}$} \put(.9,.5){$_{x_2}$}
\put(1.8,.5){$_{x_3}$}
\put(1.4,.7){$_{\az_2}$}
\put(2.4,.7){$_{\az_3}$}
\put(2.4,-.7){$_{\az_6}$}
\put(.4,.35){$_{\az_1}$}
\put(.4,-.35){$_{\az_n}$}
\psline{-}(.5,0.1)(.9,.4)
\psline{-}(.5,-.1)(.9,-.4)
\psline{-}(1.2,.5)(1.8,.5)
\psline{-}(2.2,.5)(2.8,.5)
\psline[linestyle=dotted,linewidth=1pt](1.2,-.5)(1.8,-.5)

\psline{-}(2.2,-.5)(2.8,-.5)
\psline{-}(3.1,.4)(3.5,.1)
\psline{-}(3.1,-.4)(3.5,-.1)
\put(3.3,.35){$_{\az_4}$}
\put(3.3,-.35){$_{\az_5}$}
\put(2.8,.5){$_{x_4}$} \put(2.8,-.5){$_{x_6}$}
\put(3.5,0){$_{x_5.}$}
\end{pspicture}
\end{center}

We recall from \cite{BR} that each string $w$ in $A$ defines a \emph{string module} $M(w)$ in $\mod A$. The underlying vector space of $M(w)$ is obtained by replacing each $x_i$
in $w$ by a copy of the field $k$.
The action of an arrow $\az$ of $Q$ on $M(w)$ is induced
by the relevant identity morphisms if $\az$ lies on $w$, and is zero otherwise.
For the zero string $0$, we let $M(0)$ be the zero module.
Each band $b$ defines a family of band modules $M(b,n,\phi)$ with $n\in \mathbb{N}$ and $\phi\in \mbox{Aut}(k^n)$  by replacing each $x_i$ in $b$ a copy of
the vector space $k^n$, and the action of an arrow $\az$ on $M(b,n,\phi)$ is induced by identity morphisms if $\az=\az_j$ for $j=1,2\ldots n-1$ and by $\phi$ if $\az=\az_n$ (see \cite{BR}).
\bigskip

Let $\ggz$ be a triangulation of the marked surface $(S,M)$, and denote by $(Q,W)$ the corresponding quiver with potential. In \cite{ABCP} the strings and bands of $J(Q,W)$ are related to the curves and simple closed curves respectively in $(S,M)$ :
For two curves $\gz',\gz$ in $(S,M)$ we denote by $I(\gz',\gz)$ the minimal intersection number of two representatives of the homotopic classes of $\gz'$ and $\gz$.
For each curve $\gz$ in $(S,M)$ with $d=\sum_{\gz'\in \ggz}I(\gz',\gz)$ we fix an orientation of $\gz$, and let $x_1, x_2,\ldots, x_d$ be the internal arcs of $\ggz$ that intersect $\gz$ in the fixed orientation of $\gz$, see Figure~1.

\begin{center}
\begin{pspicture}(-3.4,-2)(12,1.5)
\put(2,-2){$_{Figure~1}$}
\pscircle[linewidth=.7pt](-2,0){1}
\put(-1.4,.8){\circle*{.1}}\put(-1.4,-.8){\circle*{.1}}\put(-2,1){\circle*{.1}}
\put(-2,-1){\circle*{.1}}\put(-1,0){\circle*{.1}}
\pscircle[linewidth=.7pt](7,0){1}
\put(6.4,.8){\circle*{.1}}\put(6.4,-.8){\circle*{.1}}\put(6,0){\circle*{.1}}
\psarc[arcsepB=4pt]{<-}(-2,0){.6}{-160}{-80}\psarc[arcsepB=4pt]{<-}(7,0){.6}{-320}{-225}
\uput[l](-1,0){$_{s(\gz)}$}\uput[r](6,0){$_{e(\gz)}$}
\put(0,1.5){\circle*{.1}}\put(5,-1.5){\circle*{.1}}\put(5,1.5){\circle*{.1}}
\put(0,-1.5){\circle*{.1}}\put(2,1.5){\circle*{.1}}\put(3,1.5){\circle*{.1}}
\put(4,1.5){\circle*{.1}}\put(2,-1.5){\circle*{.1}}\put(3,-1.5){\circle*{.1}}
\psline[linewidth=.5pt]{-}(-1,0)(0,1.5)(2,1.5)(3,1.5)(2,-1.5)(0,-1.5)(-1,0)
\psline[linewidth=.5pt]{-}(6,0)(5,1.5)(4,1.5)(3,-1.5)(5,-1.5)(6,0)
\psline[linewidth=.5pt]{-}(5,1.5)(3,-1.5)(4,1.5)\psline[linewidth=.5pt]{-}(0,-1.5)(0,1.5)(2,-1.5)(2,1.5)
\psline[linewidth=.5pt]{-}(5,-1.5)(5,1.5)\psline[linewidth=.5pt]{-}(2,-1.5)(3,1.5)
\psline[linestyle=dotted,linewidth=1pt](2,-1.5)(3,-1.5)
\psline[linestyle=dotted,linewidth=1pt](3,1.5)(4,1.5)

\psline[linecolor=red](-1,0)(6,0)
{\red\uput[u](3,0){$_\gz$}}
\psline[linestyle=dotted,linewidth=1pt](3.5,0.5)(2.8,0.5)
\uput[l](0.2,0.2){$_{x_1}$}\uput[r](4.9,0.2){$_{x_d}$}
\uput[r](.8,0.2){$_{x_2}$}\uput[r](4.1,0.2){$_{x_{d-1}}$}
\uput[l](2.2,0.2){$_{x_3}$}\uput[l](-.5,-0.7){$_{x_0}$}
\put(.4,-0.7){$_{\vartriangle_1}$}\put(-.4,-0.5){$_{\vartriangle_0}$}
\put(1.3,0.9){$_{\vartriangle_2}$}\put(4.2,-0.7){$_{\vartriangle_{d-1}}$}
\put(2.3,0.9){$_{\vartriangle_3}$}\put(5.2,-0.5){$_{\vartriangle_{d}}$}
\end{pspicture}
\end{center}

Here we denote by $s(\gz) = \gz(0)$ the start point of $\gz$, and by $e(\gz)=\gz(1)$ its endpoint. Start and end point of $\gz$ lie on the boundary, indicated by the circles in Figure 1.
Along its way, the curve $\gz$ is passing through (not necessarily distinct) triangles $ \vartriangle_0, \vartriangle_1, \ldots, \vartriangle_d$. Thus we obtain a string $w(\gz)$ in $J(Q,W)$:

$$w(\gz) : x_1 \frac{~~~\az_1~}{}x_2 \frac{~~~\az_2~}{}\cdots x_{d-2} \frac{\az_{d-2}}{}x_{d-1}\frac{~~~\az_{d-1}~}{}x_{d}.$$
We recall from \cite{ABCP} the following result concerning the map $\gz \mapsto w(\gz)$:

\begin{thm}\label{ABCP} Let $\ggz$ be a triangulation of a marked surface $(S,M)$.
There exists a bijection $\gz \mapsto w(\gz)$ between the homotopy
classes of non-contractible curves in $(S,M)$ not homotopic to an arc in $\ggz$ and the strings of J(Q,W).
\end{thm}

Similarly, each closed curve $b$  in $S$ defines a cyclic walk in $J(Q,W)$.  If $b$ is not a proper power of any element in $\Pi_1(S,M)$, then it defines a band $w(b)$.

\begin{example}\label{ex1}We consider an annulus with two marked points on each boundary and a triangulation $\ggz$ with internal arcs $1,2,3,4,5$ as follows:

\begin{center}
\begin{pspicture}(-4,-2.3)(4,2)
\put(-.7,-2.3){$_{Figure~2}$}
\put(0,.5){\circle*{.1}}\put(0,-.5){\circle*{.1}}
\pscircle[linewidth=.7pt](0,0){.5}
\pscircle[linewidth=.7pt](0,0){1.8}
\put(0,1.8){\circle*{.1}}
\put(0,-1.8){\circle*{.1}}
\put(-1.8,0){\circle*{.1}}
\pscurve[linewidth=.7pt]{-}(0,1.8)(-1.2,0)(0,-1.8)
\uput[l](-1.2,0){$_{5}$}
\psarc[arcsepB=4pt]{<-}(0,0){.3}{-200}{-80}\psarc[arcsepB=4pt]{->}(0,0){2}{-200}{-175}
\psline[linewidth=.7pt]{-}(0,0.5)(0,1.8)\psline[linewidth=.7pt]{-}(0,-0.5)(0,-1.8)
\pscurve[linewidth=.7pt]{-}(0,1.8)(.8,0.5)(.75,-.2)(0,-.5)\pscurve{-}(0,1.8)(-.8,0.4)(-.75,-.2)(0,-.5)
\pscurve[linewidth=.5pt,linecolor=red]{-}(0,1.8)(1.2,0)(0,-1.4)(-1,0)(0,.8)(.6,0)(0,-.8)(-.6,0)(-.4,.4)(0,.5)
\uput[r](.7,.5){$_1$}
\uput[l](0,1.2){$_4$}
\uput[r](-.1,1.2){$_2$}
\uput[r](-.1,-1){$_3$}{\red\uput[r](.5,-.8){$_\gz$}}
\end{pspicture}
\end{center}
The associated quiver $Q_\ggz$ is cluster-tilted of type $\tilde{A}_4$ with potential $W_\ggz=\az\bz\theta:$

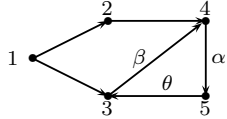
\begin{center}
\begin{pspicture}(-4,-.7)(4,.7)
\put(0,.5){\circle*{.1}}\put(0,-.5){\circle*{.1}}
\put(1.3,0.5){\circle*{.1}}\put(1.3,-0.5){\circle*{.1}}
\put(-1,0){\circle*{.1}}
\uput[l](-1,0){$_1$}\uput[u](1.3,0.4){$_4$}\uput[u](0,0.4){$_2$}\uput[d](0,-.4){$_3$}\uput[d](1.3,-0.4){$_5$}

\psline[linewidth=.7pt]{->}(-1,0)(0,.5)
\psline[linewidth=.7pt]{->}(-1,0)(0,-.5)
\psline[linewidth=.7pt]{->}(0,.5)(1.3,0.5)
\psline[linewidth=.7pt]{<-}(0,-.5)(1.3,-0.5)\psline[linewidth=.7pt]{->}(1.3,.5)(1.3,-0.5)
\psline[linewidth=.7pt]{->}(0,-.5)(1.25,0.47)
\uput[r](1.2,0){$_\az$}\uput[u](0.8,-.6){$_\theta$}
\uput[l](.7,0){$_\bz$}
\end{pspicture}
\end{center}

The fundamental group of an annulus is isomorphic to $\mathbb{Z}$, and we fix a simple closed curve $b$ as generator. From Figure~2 we obtain $w(\gz)$ and $w(b)$ as follows:
$$w(\gz): ~~3\lra 4 \longleftarrow 2 \longleftarrow 1 \lra 3 \lra 4,$$

\medskip

\begin{center}
\begin{pspicture}(-3.4,-.5)(6,0.5)
\put(-1,0){$w(b):$}
\put(.3,0){$_1$}
\psline{->}(.5,0.1)(.9,.4)
\psline{->}(.5,-.1)(.9,-.4)
\psline{->}(1.2,.4)(1.6,.1)
\psline{->}(1.2,-.4)(1.6,-.1)
\put(.95,.5){$_2$} \put(.95,-.5){$_3$}
\put(1.6,0){$_4.$}
\end{pspicture}
\end{center}
\end{example}

{\em Proof of Theorem \ref{indec}.}
Let $\ggz$ be a triangulation of $(S,M)$.
As explained in section \ref{quiver},  the  indecomposable objects in $\mc_{(S,M)}$ are either given by indecomposable modules over the Jacobian algebra $J(Q,W)$, or they correspond to the indecomposable summands of $T_\ggz$, thus to the internal arcs in $\Gamma$.
The (finite-dimensional) indecomposable modules over a string algebra $A$ are classified in \cite{BR}: Each indecomposable $A-$module is (isomorphic to) a string or a band module. The string module $M(w)$ is isomorphic to the string module $M(w^{-1})$ defined by the inverse string $w^{-1}$, and the band module $M(b,n,\phi)$ is isomorphic to $M(b',n,\phi)$ whenever $b'$ is obtained from $b$ by inversing or cyclic permutation. Apart from that, there are no isomorphisms between string or band modules.
\medskip

Thus each non-contractible curve $\gz$ which is not homotopic to a boundary segment of $(S,M)$ corresponds to an indecomposable object in  $\mc_{(S,M)}$. In case $\gz$ is not an internal arc in $\Gamma$, it corresponds to the string module $M(w(\gz))$.
From Theorem \ref{ABCP} and what we have discussed above
we conclude that two such curves $\gz,\delta$ are isomorphic as objects in $\mc_{(S,M)}$ precisely when $\gz$ is homotopic to $\delta$ or to its inverse $\delta^{-1}$. We refer to these objects as the string objects or curves in $(S,M),$ as described in part (1) of Theorem \ref{indec}.
\medskip

The remaining indecomposable objects in $\mc_{(S,M)}$ correspond to the band modules $M(b,n,\phi)$ over $J(Q,W)$, we refer to them as band objects.
They are parametrized by a positive integer $n$, an automorphism $\phi$ of $k^n$ which is given by an element of $k^*$ since $k$ is algebraically closed, and a band $b$ of $J(Q,W)$.
The fundamental group $\Pi_1(S,M)$ is a free group with a finite number of generators which are given by closed simple curves in $S$.
In order to avoid counting curves with opposite orientation twice we consider the elements in  $\Pi_1(S,M)$ up to the equivalence relation $a \sim a^{-1}$.
\medskip

Moreover, to comply with the definition of a band module, we write each element $a$  in $\Pi_1(S,M)$ as $a=b^n$ (multiplicatively written) for some $b \in \Pi_1(S,M)$ which itself is not a proper power of an element in $\Pi_1(S,M)$. Furthermore, we consider the elements of $\Pi_1(S,M)$ up to cyclic permutation of their factors.
Then it is clear that the band modules $M(b,n,\phi)$ over $J(Q,W)$ correspond bijectively to  $k^* \times {\Pi^*_1 (S,M)}/ \sim$, where  $\Pi^*_1 (S,M)/\sim$  is given by the non-zero elements of the fundamental group of $(S,M)$ subject to the equivalence relation generated by $a \sim a^{-1}$ and cyclic permutation.
This is the description of the band objects in $\mc_{(S,M)}$ given in part (2) of Theorem \ref{indec}.
\hfill $\Box$

\section{Irreducible morphisms in $\mc_{(S,M)}$}
Based on the geometric characterization of the indecomposable objects in $\mc_{(S,M)}$ in the previous section, we study in this section the irreducible morphisms in $\mc_{(S,M)}$.
\subsection{The AR-quiver of a string algebra.}
We first recall some basic definitions from \cite{BR}. Let $A=kQ/I$ be a finite-dimensional string algebra with $Q=(Q_0,Q_1)$ and $\mathcal{S}$ the set of all strings in $A$. A string $w$ \emph{starts \emph{(}or ends\emph{)}  on a peak} if there is no arrow $\az \in Q_1$ with $w\az \in \ms$ (or $ \az^{-1}w \in \ms$); likewise, a string $w$ \emph{starts \emph{(}or ends\emph{)} in a deap} if there is no arrow $\bz \in Q_1$ with $w\bz^{-1} \in \ms$ (or $ \bz w \in \ms$).

A string $w=\az_1\az_2\cdots\az_n$ with all $\az_i \in Q_1$ is called \emph{direct string}, and a string of the form $w^{-1}$ where $w$ is a direct string is called \emph{inverse string}.
Strings of length zero are both direct and inverse. For each arrow $\az \in Q_1$, let $N_\az=U_\az \az V_\az$ be the unique string such that $U_\az$ and $V_\az$ are inverse strings and $N_\az$ starts in a deep and ends on a peak (see the following figure).

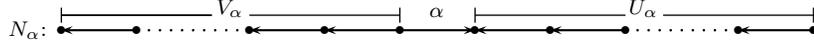
\begin{center}
\begin{pspicture}(1,-.1)(10,.5)
\put(0,0){\circle*{.1}}
\put(1,0){\circle*{.1}}\put(9,0){\circle*{.1}}\put(10,0){\circle*{.1}}
\put(2.5,0){\circle*{.1}}\put(3.5,0){\circle*{.1}}\put(4.5,0){\circle*{.1}}
\put(6.5,0){\circle*{.1}}\put(5.5,0){\circle*{.1}}\put(7.5,0){\circle*{.1}}
\psline{<-}(0,0)(1,0)\psline{<-}(2.5,0)(3.5,0)\psline{<-}(3.5,0)(4.5,0)
\psline{<-}(9,0)(10,0)\psline{<-}(6.5,0)(7.5,0)\psline{<-}(5.5,0)(6.5,0)
\psline{->}(4.5,0)(5.5,0)
\psline[linestyle=dotted,linewidth=1pt](1,0)(2.5,0)\psline[linestyle=dotted,linewidth=1pt](9,0)(7.5,0)

\psline[linewidth=.5pt](0,.2)(2,.2)\psline[linewidth=.5pt](2.5,.2)(4.5,.2)
\psline[linewidth=.5pt](5.5,.2)(7.5,.2)
\psline[linewidth=.5pt](8,.2)(10,.2)\psline[linewidth=.5pt](0,.2)(2,.2)
\psline[linewidth=.5pt](0,.1)(0,.3)\psline[linewidth=.5pt](4.5,.1)(4.5,.3)
\psline[linewidth=.5pt](5.5,.1)(5.5,.3)\psline[linewidth=.5pt](10,.1)(10,.3)
\uput[u](5,0){$_{\az}$}
\uput[u](2.25,0){$_{V_{\az}}$}\uput[u](7.75,0){$_{U_{\az}}$}
\uput[l](0,0){$_{N_{\az}:}$}
\end{pspicture}
\end{center}
If the string $w$ does not start on a peak, we define $w_h=w\az V_\az$ and say that $w_h$ is obtained from $w$ by adding a hook on the starting point $s(w)$. Dually, if $w$ does not end on a peak , we define $_hw={V_\az}^{-1}\az^{-1}w$ and say that $_hw$ is obtained from $w$ by \emph{adding a hook }on the ending point $e(w)$ (see the following figure).

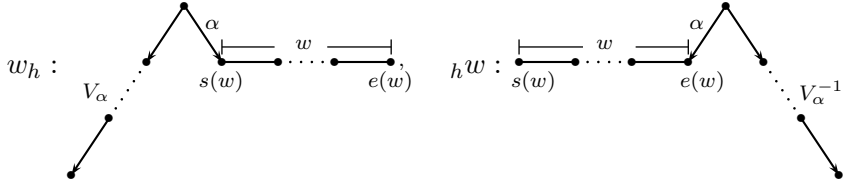
\begin{center}
\begin{pspicture}(-6.2,-1.7)(6,.85)
\put(-.95,0){\circle*{.1}}\put(-.85,0){$_{,}$}\put(-1.7,0){\circle*{.1}}\put(-2.45,0){\circle*{.1}}
\put(-3.2,0){\circle*{.1}}\put(-4.2,0){\circle*{.1}}\uput[d](-3.2,.05){$_{s(w)}$}\uput[d](3.2,0.05){$_{e(w)}$}
\uput[d](-.95,0.05){$_{e(w)}$}\uput[d](.95,0.05){$_{s(w)}$}
\put(-3.7,0.75){\circle*{.1}}\put(-4.7,-.75){\circle*{.1}}
\put(-5.2,-1.5){\circle*{.1}}\uput[l](-5.2,-.1){${w_h:}$}
\psline{-}(-.95,0)(-1.7,0)
\psline[linestyle=dotted,linewidth=1pt](-1.7,0)(-2.45,0)\psline{-}(-2.45,0)(-3.2,0)
\psline{->}(-3.7,0.75)(-3.2,0)\psline{->}(-3.7,0.75)(-4.2,0)\psline[linestyle=dotted,linewidth=1pt](-4.2,0)(-4.7,-0.75)
\psline{->}(-4.7,-.75)(-5.2,-1.5)
\psline[linewidth=.5pt](-0.95,.2)(-1.8,.2)\psline[linewidth=.5pt](-2.5,.2)(-3.2,.2)
\psline[linewidth=.5pt](-0.95,.1)(-0.95,.3)\psline[linewidth=.5pt](-3.2,.1)(-3.2,.3)
\uput[r](-3.6,0.5){$_{\az}$}
\uput[l](-4.5,-.4){$_{V_{\az}}$}
\uput[u](-2.1,0){$_w$}
\put(.75,0){\circle*{.1}}\put(1.5,0){\circle*{.1}}\put(2.25,0){\circle*{.1}}
\put(3,0){\circle*{.1}}
\put(4,0){\circle*{.1}}
\put(3.5,0.75){\circle*{.1}}
\put(4.5,-.75){\circle*{.1}}\put(5,-1.5){\circle*{.1}}
\uput[l](0.7,-.1){${_h{w:}}$}
\psline{-}(.75,0)(1.5,0)
\psline[linestyle=dotted,linewidth=1pt](1.5,0)(2.25,0)\psline{-}(2.25,0)(3,0)
\psline{->}(3.5,0.75)(3,0)\psline{->}(3.5,0.75)(4,0)\psline[linestyle=dotted,linewidth=1pt](4,0)(4.5,-0.75)
\psline{->}(4.5,-.75)(5,-1.5)
\psline[linewidth=.5pt](0.75,.2)(1.6,.2)\psline[linewidth=.5pt](2.3,.2)(3,.2)
\psline[linewidth=.5pt](0.75,.1)(0.75,.3)\psline[linewidth=.5pt](3,.1)(3,.3)
\uput[l](3.4,0.5){$_{\az}$}
\uput[r](4.3,-.4){$_{V^{-1}_{\az}}$}
\uput[u](1.9,0){$_{w}$}
\end{pspicture}
\end{center}

Suppose now that the string $w$ starts on a peak. If $w$ is not a direct string, we can write $w=w_c\bz^{-1}\gz_1\gz_2\cdots\gz_r=w_c\bz^{-1} U_\bz^{-1}$ for some $\bz\in Q_1$ and $r\geq 0$. We say in this case that $w_c$ is obtained from $w$ by \emph{deleting a cohook} on $s(w)$. If $w$ is a direct string, we define $w_c=0$.
Dually, assume that $w$ ends on a peak. Then, if $w$ is not an inverse string, we can write $w=\gz_r^{-1}\cdots\gz_2^{-1}\gz_1^{-1}\bz_cw=U_\bz\bz_cw$ for some $\bz\in Q_1$ and $r\geq 0$. We say that $_cw$ is obtained from $w$ by deleting a cohook on $e(w)$, and if $w$ is an inverse string, we define $_cw=0$ (see the following figure).

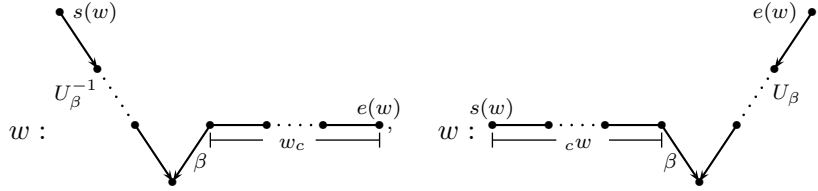
\begin{center}
\begin{pspicture}(-6,-.85)(6,1.7)
\uput[r](-5,1.5){$_{s(w)}$}\uput[l](5,1.5){$_{e(w)}$}\uput[u](-.75,-.1){$_{e(w)}$}\uput[u](.75,-.1){$_{s(w)}$}
\put(-.75,0){\circle*{.1}}\put(-1.5,0){\circle*{.1}}\put(-2.25,0){\circle*{.1}}
\put(-3,0){\circle*{.1}}
\put(-4,0){\circle*{.1}}
\put(-3.5,-0.75){\circle*{.1}}
\put(-4.5,.75){\circle*{.1}}\put(-5,1.5){\circle*{.1}}
\uput[l](-5,-.1){${w:}$}
\psline{-}(-.75,0)(-1.5,0)
\psline[linestyle=dotted,linewidth=1pt](-1.5,0)(-2.25,0)\psline{-}(-2.25,0)(-3,0)
\psline{<-}(-3.5,-0.75)(-3,0)\psline{<-}(-3.5,-0.75)(-4,0)\psline[linestyle=dotted,linewidth=1pt](-4,0)(-4.5,0.75)
\psline{<-}(-4.5,.75)(-5,1.5)
\psline[linewidth=.5pt](-0.75,-.2)(-1.6,-.2)\psline[linewidth=.5pt](-2.3,-.2)(-3,-.2)
\psline[linewidth=.5pt](-0.75,-.1)(-0.75,-.3)\psline[linewidth=.5pt](-3,-.1)(-3,-.3)
\uput[r](-3.4,-0.5){$_{\bz}$}
\uput[l](-4.3,.4){$_{U^{-1}_{\bz}}$}
\uput[d](-1.9,0){$_{w_c}$}\put(-.65,0){$_{,}$}

\put(.75,0){\circle*{.1}}\put(1.5,0){\circle*{.1}}\put(2.25,0){\circle*{.1}}
\put(3,0){\circle*{.1}}\put(4,0){\circle*{.1}}\put(3.5,-0.75){\circle*{.1}}
\put(4.5,.75){\circle*{.1}}\put(5,1.5){\circle*{.1}}\uput[l](0.7,-.1){${w:}$}
\psline{-}(.75,0)(1.5,0)
\psline[linestyle=dotted,linewidth=1pt](1.5,0)(2.25,0)\psline{-}(2.25,0)(3,0)
\psline{<-}(3.5,-0.75)(3,0)\psline{<-}(3.5,-0.75)(4,0)\psline[linestyle=dotted,linewidth=1pt](4,0)(4.5,0.75)
\psline{<-}(4.5,.75)(5,1.5)
\psline[linewidth=.5pt](0.75,-.2)(1.6,-.2)\psline[linewidth=.5pt](2.3,-.2)(3,-.2)
\psline[linewidth=.5pt](0.75,-.1)(0.75,-.3)\psline[linewidth=.5pt](3,-.1)(3,-.3)
\uput[l](3.4,-0.5){$_{\bz}$}
\uput[r](4.3,.4){$_{U_{\bz}}$}
\uput[d](1.9,0){$_{_cw}$}
\end{pspicture}
\end{center}

The following theorem summarizes the results from  \cite{BR} concerning AR-sequences between string modules:

\begin{thm}[\cite{BR}]\label{BR} For  a string algebra $A$, let $w$ be a string such that  $M(w)$ is not an injective $A-$module. Then the AR-sequence starting in $M(w)$ is given
\begin{itemize}
\item[(1)] if $w$ neither starts nor ends on a peak by
 $$0\lra M(w) \lra M(w_h)\oplus M(_hw) \lra M(_hw_h) \lra 0,$$

\item[(2)] if $w$ does not start but ends on a peak by
$$0\lra M(w) \lra M(w_h)\oplus M(_cw) \lra M(_cw_h) \lra 0,$$

\item[(3)] if $w$ starts but does not end on a peak by
$$0\lra M(w) \lra M(w_c)\oplus M(_hw) \lra M(_hw_c) \lra 0,$$

\item[(4)] if $w$ both starts and ends on a peak by
$$0\lra M(w) \lra M(w_c)\oplus M(_cw) \lra M(_cw_c) \lra 0.$$
\end{itemize}
\end{thm}

\begin{rem}\label{rmk1}
\begin{itemize}
\item[(1)]Each arrow $\az \in Q_1$ defines a string $N_\az=U_\az \az V_\az$ which starts in a deep and ends on a peak. Then $w=U_\az$ does not start but ends on a peak, hence we have an AR-sequence by the above theorem:
    $$0\lra M(w) \lra M(w_h)\oplus M(_cw) \lra M(_cw_h) \lra 0.$$
Here $M(_cw)=0$ since $w$ is an inverse string, and $w_h=N_\az$ and $_cw_h={_c{(N_\az)}}=V_\az$ by definition. Hence, for each $\az\in Q_1$, there is an AR-sequence with an indecomposable middle term in $\mod A$:
$$0 \lra M(U_\az) \lra M(N_\az)\lra M(V_\az)\lra 0.$$
In fact, the AR-sequences between string modules which admit only one direct summand in the middle term are indexed by the arrows in $Q$ (see \cite{BR}).
\bigskip

\item[(2)]If the string module $M(w)$ is injective, then $w$ both starts and ends on a peak, thus we can write
$$w=\gz_r^{-1}\gz_{r-1}^{-1}\cdots\gz_1^{-1}\bz_1\bz_2\cdots\bz_{s-1}\bz_s \mbox { or}$$
$$w=U_{\bz_1}\bz_1\bz_2\cdots\bz_{s-1}\bz_s=\gz_r^{-1}\gz_{r-1}^{-1}\cdots
    \gz_1^{-1}U_{\gz_1}^{-1}$$ where $r+s\geq 1$, hence $_cw=\gz_r^{-1}\gz_{r-1}^{-1}\cdots\gz_2^{-1}, w_c=\bz_2\cdots\bz_{s-1}\bz_s$, and $M(w)\diagup{Soc(M(w))}=M(_cw)\oplus M(w_c)$ which yields two irreducible morphisms $M(w)\lra M(_cw),M(w)\lra M(w_c)$. Moreover, $_cw_c={_c(w_c)}=(_cw)_c$ is the empty string, hence $M(_cw_c)$ is the zero module. Thus one might say  that the case where M(w) is injective is contained in the case (4) of the above theorem.
    \bigskip

\item[(3)] The number of indecomposable summands in the middle term of an AR-sequence between string modules is  at most two.
\end{itemize}
\end{rem}
\medskip

\subsection{Irreducible morphisms in $J(Q,W)$}
We now fix  a triangulation $\ggz=\{\tau_1,\tau_2,\ldots \tau_n, \tau_{n+1},\ldots,\tau_{n+m}\}$ of the marked surface $(S,M)$ where $\tau_1,\ldots,\tau_n$ are internal arcs and $\tau_{n+1},\ldots,\tau_{n+m}$ are boundary arcs.
The aim of this subsection is to describe the AR-quiver of the string algebra $J(Q,W)$ in terms of objects in  $\mc_{(S,M)}$, which we identified with curves and closed curves in $(S,M)$.
A curve $\gz$ in $(S,M)$ defines a string $w(\gz)$ in $J(Q,W)$ (which is empty if $\gz$ is contractible or homotopic to an arc in $\ggz$), and consequently a string module $M(w(\gz))$ in $\mod J(Q,W)$ (which is zero if $\gz$ is contractible or homotopic to an arc in $\ggz$). We use the notation $M(\gz)$ for $M(w(\gz))$ and from the discussion in Section 2 we know that the string module  $M(\gz)$ is the image of $\gz$ under the projection functor
$\Ext_{\mc_\ggz}^1(T_\ggz,-): \mc_{(S,M)} \to \mod J(Q,W)$.

To define elementary moves on curves in $(S,M)$ we use the fact that the orientation of $S$ induces an orientation on each boundary component of $S$:
For any curve $\gz$ in $(S,M)$ we denote by $_s\gz$ the \emph{Pivot elementary move} of $\gz$ on its starting point, meaning that the curve $_s\gz$ is obtained from $\gz$ by moving the starting point $s(\gz)$
clockwise to the next marked point $b$ on the same boundary (note that $b=s(\gz)$ if there is only one marked point lying on the same boundary). Similarly, we denote by $\gz_e$ the Pivot elementary move of $\gz$ on its ending point, see Figure~3. Iterated Pivot elementary moves are denoted  $_s\gz_e={_s(\gz_e)}=(_s\gz)_e$, $ {_{s^2}\gz}={_s{(_s\gz)}}$ and $\gz_{e^2}={(\gz_e)}_e$, respectively.

\begin{center}
\begin{pspicture}(-3.4,-1.2)(7,1.3)
\put(.9,-1.4){$_{Figure~3}$}
\pscircle[linewidth=.7pt](-.5,0){1}
\put(.1,.8){\circle*{.1}}
\put(.1,.8){\circle*{.1}}\put(.1,-.8){\circle*{.1}}
\put(-.5,1){\circle*{.1}}\put(-.5,-1){\circle*{.1}}
\pscircle[linewidth=.7pt](4,0){1}
\put(3.4,.8){\circle*{.1}}
\put(3.4,-.8){\circle*{.1}}
\put(5,0){\circle*{.1}}
\psline[linestyle=dotted,linewidth=1pt](.1,-.8)(3.4,.8)
\psline[linecolor=red](.1,.8)(3.4,-.8)
\psline[linewidth=.5pt](.1,-.8)(3.4,-.8)
\psline[linewidth=.5pt](.1,.8)(3.4,.8)
\uput[u](1.6,.7){$_{\gz_e}$}
\uput[d](1.6,-.7){$_{_s\gz}$}
{\red\uput[u](1.4,0.1){$_{\gz}$}}
\uput[u](1.1,-0.4){$_{_s\gz_e}$}
\uput[l](.2,0.6){$_{s(\gz)}$}
\uput[r](3.5,-.7){$_{e(\gz)}$}
\psarc[arcsepB=4pt]{<-}(-.5,0){.6}{-160}{-80}\psarc[arcsepB=4pt]{<-}(4,0){.6}{-320}{-225}
\uput[l](.1,-.7){$_b$}
\end{pspicture}
\end{center}

\begin{lem}\label{lem1} Let $\ggz$ be a triangulation of a marked surface $(S,M)$, and let $\gz$ be a curve in $(S,M)$ such that the string $w(\gz)$ of $J(Q,W)$ is non-empty.
\begin{itemize}
\item[(1)]If $w(\gz)$ does not start on a peak, then $w(_s\gz)$ is obtained by adding a hook on $s(w(\gz))$. Thus $w(_s\gz)=w(\gz)_h$  and there is an irreducible morphism in $\mod J(Q,W): M(\gz)\lra M(_s\gz).$
\item[(2)]If $w(\gz)$ does not end on a peak, then $w(\gz_e)$ is obtained by adding a hook on $e(w(\gz))$. Thus $w(\gz_e)={_hw(\gz)}$ and there is an irreducible morphism in $\mod J(Q,W): M(\gz)\lra M(\gz_e).$
\end{itemize}
\end{lem}
\begin{pf} We only prove $(1)$ since $(2)$ is obtained dually.

\begin{center}
\begin{pspicture}(-3.4,-2.8)(12,1.3)
\put(2,-2.8){$_{Figure~4}$}
\pscircle[linewidth=.7pt](-2,0){1}
\put(-1.4,.8){\circle*{.1}}\put(-1.4,-.8){\circle*{.1}}
\put(-2,1){\circle*{.1}}\put(-2,-1){\circle*{.1}}
\put(-1,0){\circle*{.1}}
\pscircle[linewidth=.7pt](7,0){1}
\put(6.4,.8){\circle*{.1}}\put(6.4,-.8){\circle*{.1}}\put(6,0){\circle*{.1}}
\psarc[arcsepB=4pt]{<-}(-2,0){.6}{-160}{-80}\psarc[arcsepB=4pt]{<-}(7,0){.6}{-320}{-225}
\uput[l](-1,0){$_{s(\gz)}$}\uput[r](6,0){$_{e(\gz)}$}
\put(0,1.5){\circle*{.1}}\put(5,-1.5){\circle*{.1}}\put(5,1.5){\circle*{.1}}
\put(2,1.5){\circle*{.1}}\put(3,1.5){\circle*{.1}}
\put(4,1.5){\circle*{.1}}\put(2,-1.5){\circle*{.1}}
\put(3,-1.5){\circle*{.1}}
\psline{-}(-1,0)(0,1.5)(2,1.5)(3,1.5)(2,-1.5)(-1,0)
\psline[linewidth=.5pt]{-}(6,0)(5,1.5)(4,1.5)(3,-1.5)(5,-1.5)(6,0)\psline[linewidth=.5pt]{-}(5,1.5)(3,-1.5)(4,1.5)
\psline[linewidth=.5pt]{-}(0,1.5)(2,-1.5)(2,1.5)\psline[linewidth=.5pt]{-}(5,-1.5)(5,1.5)
\psline[linewidth=.5pt]{-}(2,-1.5)(3,1.5)\psline[linestyle=dotted,linewidth=1pt](2.3,-1.5)(3,-1.5)
\psline[linestyle=dotted,linewidth=1pt](3,1.5)(4,1.5)
\psline[linecolor=red](-1,0)(6,0)
{\red\uput[u](3,0){$_\gz$}}
\psline[linestyle=dotted,linewidth=1pt](3.5,0.5)(2.8,0.5)
\uput[r](4.9,0.2){$_{\tau_{i_d}}$}\uput[r](.8,0.2){$_{\tau_{i_1}}$}
\uput[r](4.1,0.2){$_{\tau_{i_{d-1}}}$}\uput[l](2.2,0.2){$_{\tau_{i_2}}$}
\uput[l](.5,-0.4){$_{\tau_{i_0}}$}\put(0,.6){$_{\vartriangle_0}$}
\put(1.3,0.9){$_{\vartriangle_1}$}\put(4.2,-0.7){$_{\vartriangle_{d-1}}$}
\put(2.3,0.9){$_{\vartriangle_2}$}\put(5.2,-0.5){$_{\vartriangle_{d}}$}
\put(-1.4,-3){\circle*{.1}}
\psline[linewidth=.5pt]{-}(-1.4,-.8)(-1.4,-3)(2,-1.5)
\psline[linewidth=.5pt]{-}(-1,0)(-1.4,-3)
\put(-.2,-2.47){\circle*{.1}}
\put(.8,-2.03){\circle*{.1}}
\put(.6,-1.6){$_{\tau_{j_1}}$}
\put(-.3,-2){$_{\tau_{j_{r-1}}}$}
\uput[l](-1.2,-1.5){$_{\tau_{j_{r+1}}}$}
\put(-1.2,-2.3){$_{\tau_{j_r}}$}
{\red\psline[linewidth=.5pt]{-}(-.2,-2.47)(-1,0)(.8,-2.03)}
\psline[linewidth=.5pt]{-}(2.3,-1.5)(2,-1.5)
\psline[linecolor=blue]{->}(.8,-.9)(1,0)
\psline[linecolor=blue]{->}(.8,-.9)(.2,-1.3)
\psline[linestyle=dotted,linecolor=blue](.2,-1.3)(-.5,-1.6)
\psline[linecolor=blue]{->}(-.5,-1.6)(-1.2,-1.8)
{\blue\uput[l](.6,-1){$_{\bz_1}$}
\uput[r](.8,-.6){$_{\az}$}
\uput[l](-.5,-1.5){$_{\bz_r}$}}
\uput[l](-1.4,-.7){$_b$}
\psline[linewidth=.4pt,linecolor=blue]{-}(-1.4,-.8)(6,0)
{\blue\uput[u](2.9,-.7){$_{_s\gz}$}}
\end{pspicture}
\end{center}
As shown in  Figure~4, we denote by $\tau_{i_1},\tau_{i_2},\ldots, \tau_{i_d}$ the internal arcs of $\ggz$ that intersect $\gz$. Let $b$ be the marked point lying clockwise next to $s(\gz)$ on the same boundary component. Since $\gz$ does not start on a peak, there exits an arrow $\az:\tau_{i_0} \to \tau_{i_1}$  in $Q$ such that  $w(\gz)\az$ is a string in $J(Q,W)$, for some internal arc $\tau_{i_0} \in \ggz$.
Let $\tau_{j_1},\tau_{j_2}, \ldots, \tau_{j_r}$ be all internal arcs in $\ggz$ which intersect $\tau_{i_0}$ in the vertex $s(\gz)$ and which are successors of $\tau_{i_0}$ with respect to clockwise orientation at the common vertex $s(\gz)$.
We denote by $\bz_1,\bz_2, \ldots, \bz_r$ the arrows of $Q$ induced by the internal arcs $\tau_{i_0},\tau_{j_1}, \ldots,\tau_{j_r}$.
Then $w(_s\gz)=w(\gz)\az\bz_1^{-1}\bz_2^{-1}\cdots\bz_r^{-1}=w(\gz)\az V_\az$ where $V_\az=\bz_1^{-1}\bz_2^{-1}\cdots\bz_r^{-1}$. This means that $w(_s\gz)$ is obtained from $w(\gz)$ by adding a hook on $s(w(\gz))$. By Theorem \ref{BR} there is an irreducible morphism: $M(\gz)\lra M(_s\gz)$ in $\mod J(Q,W).$
\end{pf}
\begin{lem}\label{lem2}  Let $\ggz$ be a triangulation of a marked surface $(S,M)$, and let $\gz$ be a curve in $(S,M)$ such that the string $w(\gz)$ of $J(Q,W)$ is non-empty.
\begin{itemize}
\item[(1)]If $w(\gz)$ starts on a peak, then $w(_s\gz)$ is obtained by deleting a cohook on $s(w(\gz))$. Thus $w(_s\gz)=w(\gz)_c$ and if $w(\gz)_c$ is non-empty there is an irreducible morphism in $\mod J(Q,W): M(\gz)\lra M(_s\gz).$
\item[(2)]If $w(\gz)$ ends on a peak, then $w(\gz_e)$ is obtained by deleting a cohook on $e(w(\gz))$. Thus  $w(\gz_e)= {_cw(\gz)}$ and if $_cw(\gz)$ is non-empty there is an irreducible morphism in $\mod J(Q,W): M(\gz)\lra M(\gz_e).$
\end{itemize}
\end{lem}
\begin{pf}As before we only prove part (1) of the Lemma.

\begin{center}
\begin{pspicture}(-3.4,-2)(12,2)
\pscircle[linewidth=.7pt](-2,0){1}
\put(-1.4,.8){\circle*{.1}}\put(-1.1,-.4){\circle*{.1}}
\put(-2,1){\circle*{.1}}\put(-2,-1){\circle*{.1}}
\put(1,1.5){\circle*{.1}}
\uput[l](-1.1,-.3){$_b$}
\pscircle[linewidth=.7pt](7,0){1}
\put(6.4,.8){\circle*{.1}}\put(6.4,-.8){\circle*{.1}}\put(6,0){\circle*{.1}}
\psarc[arcsepB=4pt]{<-}(-2,0){.6}{-160}{-80}\psarc[arcsepB=4pt]{<-}(7,0){.6}{-320}{-225}
\uput[l](-1.2,.6){$_{s(\gz)}$}\uput[r](6,0){$_{e(\gz)}$}
\put(0,1.5){\circle*{.1}}\put(5,-1.5){\circle*{.1}}\put(5,1.5){\circle*{.1}}
\put(0,-1.5){\circle*{.1}}\put(2,1.5){\circle*{.1}}\put(3.5,1.5){\circle*{.1}}
\put(4,1.5){\circle*{.1}}\put(2,-1.5){\circle*{.1}}\put(3,-1.5){\circle*{.1}}
\psline[linewidth=.5pt]{-}(-1.4,.8)(0,1.5)(2,1.5)(3.5,1.5)(2,-1.5)(0,-1.5)
\psline[linewidth=.5pt]{-}(6,0)(5,1.5)(4,1.5)(3,-1.5)(5,-1.5)(6,0)
\psline[linewidth=.5pt]{-}(5,1.5)(3,-1.5)(4,1.5)
\psline[linewidth=.5pt]{-}(-1.1,-.4)(1,1.5)
\psline[linewidth=.5pt]{-}(2,1.5)(-1.1,-.4)(0,1.5)
\psline[linewidth=.5pt]{-}(0,-1.5)(-1.1,-.4)(3.5,1.5)
\psline[linewidth=.5pt]{-}(5,-1.5)(5,1.5)
\psline[linewidth=.5pt]{-}(0,-1.5)(3.5,1.5)
\psline[linestyle=dotted,linewidth=1pt](2,-1.5)(3,-1.5)
\psline[linestyle=dotted,linewidth=1pt](3.5,1.5)(4,1.5)
\psline[linewidth=1pt,linestyle=dotted]{-}(2.9,0)(3.3,0)

{\red\uput[u](3.2,.2){$_\gz$}}{\blue\uput[u](3,-.6){$_{_s\gz}$}}
\uput[r](4.9,0.4){$_{\tau_{i_d}}$}
\uput[r](5.3,0.8){$_{\tau_{i_{d+1}}}$}
\uput[r](3.7,-0.5){$_{\tau_{i_{d-1}}}$}

\psline[linewidth=.6pt,linecolor=red]{-}(-1.4,.8)(6,0)
\psline[linewidth=.6pt,linecolor=blue]{-}(-1.1,-.4)(6,0)
\psline[linewidth=.8pt,linecolor=red]{->}(-.2,1)(.4,1)
\psline[linestyle=dotted,linewidth=.8pt,linecolor=red](.5,1)(1.1,1)
\psline[linewidth=.8pt,linecolor=red]{->}(1.1,.9)(1.2,.6)
\psline[linewidth=.8pt,linecolor=red]{->}(1.4,-.3)(1.2,.6)
{\red\put(1.4,0.2){$_{\az}$}
\put(1.2,0.8){$_{\bz_r}$}
\put(.1,1.2){$_{\bz_1}$}}
\uput[l](-.9,.3){$_{\tau_{i_0}}$}
\uput[l](-.4,.5){$_{\tau_{i_1}}$}
\uput[d](-.4,-.5){$_{\tau_{i_{-1}}}$}
\uput[l](0.1,.5){$_{\tau_{i_2}}$}
\uput[r](.2,.1){$_{\tau_{i_{r+1}}}$}
\uput[d](1.6,-.3){$_{\tau_{i_{r+2}}}$}
\psline[linestyle=dotted]{-}(0,.5)(.4,0.5)
\put(2,-2){$_{Figure~5}$}
\end{pspicture}
\end{center}

We denote by $b$ the marked point lying clockwise next to $s(\gz)$ on the same boundary.
We further denote by $\tau_{i_1},\tau_{i_2}, \ldots, \tau_{i_d}$  the internal arcs that intersect $\gz$ in the order indicated in Figure~5.
Since $w(\gz)$ starts on a peak, there is an $r \ge 1$ such that the arcs $\tau_{i_1},\tau_{i_2}, \ldots, \tau_{i_{r+1}}$ intersect in the vertex $b$ and induce arrows $\bz_1,\bz_2, \ldots, \bz_{r}$ in $Q_1$ as shown in Figure~5.
We choose $r$ to be maximal and distinguish the following cases:
\begin{itemize}
\item[(1)]If $r+1=d$ then $w(\gz)$ is a direct string and $w({_s\gz})=w(\gz)_c$ is the empty string.
\item[(2)]If $r+1 < d$ then by maximality of $r$ there is an arrow $\alpha$  in $Q_1$ from $\tau_{i_{r+2}}$ to $\tau_{i_{r+1}}$.
Since $\gz$ starts on a deep, $\tau_{i_0}$ is a boundary arc which implies $U_\az^{-1}=\bz_r\bz_{r-1}$ $\cdots\bz_1$. Moreover, we know $w(\gz)=w_0\az^{-1}\bz_r\bz_{r-1}\cdots\bz_1=w_0\az^{-1}U_\az^{-1}$.
Thus $w({_s\gz})=w_0$  which means that $w(_s\gz)$ is obtained by deleting a cohook on $s(w(\gz))$. If $w(_s\gz)$ is non-empty, we obtain by Theorem \ref{BR} that  there is an irreducible morphism in $\mod J(Q,W): M(\gz)\lra M(_s\gz).$
\end{itemize}
\end{pf}

\begin{thm}\label{thm1}Let $\ggz$ be a triangulation of a marked surface $(S,M)$, and let $\gz$ be a curve in $(S,M)$.
Then each irreducible morphism in $\mod J(Q,W)$ starting in $M(\gz)$ is obtained by Pivot elementary moves on endpoints of $\gz$. Moreover, all AR-sequences between string modules in $\mod J(Q,W)$ are of the form
$$0\lra M(\gz)\lra M(_s\gz)\oplus M(\gz_e)\lra  M(_s\gz_e)\lra 0$$
for some curve $\gz$ in $(S,M)$.
\end{thm}
\begin{pf} The irreducible morphisms and AR-sequences between string modules are described in Theorem \ref{BR}.
Among the four cases listed there, we only consider the second case here, the others being similar:
Suppose $w(\gz)$ does not start but ends on a peak, then we get from Lemma \ref{lem1} and Lemma \ref{lem2} two irreducible morphisms
$M(\gz)\lra M(_s\gz)$ and $M(\gz)\lra M(\gz_e).$
Moreover, the construction of $_s\gz$ does not change any information of $\gz$ on $e(\gz)$, that is $w(_s\gz)$ also ends on a peak if it is non-empty.
By Lemma \ref{lem2} again we get an irreducible morphism $M(_s\gz)\lra M(_s\gz_e)$ in case $w(_s\gz_e)$ is non-empty, and similarly, there is an irreducible morphism $M(\gz_e)\lra M(_s\gz_e)$ induced by $\gz_e$. This completes the proof since the number of summands in the middle term of an AR-sequence between string modules is at most two.
\end{pf}
To complete the description of irreducible morphisms in $\mod J(Q,W)$, we recall from \cite{BR} that  each band in $J(Q,W)$ yields a $k^*-$family of homogeneous tubes in the AR-quiver of $J(Q,W)$, given by an embedding
 $\mod k[t,t^{-1}] \to \mod J(Q,W)$.

\begin{cor}\label{ar-translation}The AR-translation in $\mod J(Q,W)$ is given by simultaneous counterclockwise rotation of starting and ending point of a curve to the next marked points at the boundary, that is
$$\tau^{-1}(M(\gz))=M(_s\gz_e)$$
if $M(\gz)$ is not an injective $J(Q,W)-$module.
\end{cor}
\begin{pf} For string modules this is shown in Theorem \ref{thm1} above. For any band $b$ in $J(Q,W)$ we have $b={_sb_e}$ as elements in $\Pi_1(S,M)$. Moreover, the corresponding band modules $M(b,n,\phi)$ lie on homogeneous tubes and thus satisfy
 $\t^{-1}(M(b,n,\phi))=M(b,n,\phi)$.
 \end{pf}

\begin{example}We reconsider the curve $\gz$ from Example \ref{ex1}. From  Figure~2 it is easy to obtain the following descriptions:

\begin{center}
\begin{pspicture}(-3.4,-.5)(6,0.5)

\put(-4.8,0){$_{M(\gz)=}$}
\put(-3.2,0.3){$_{~1}$}
\put(-3.5,0){$_{3~2~3~,}$}
\put(-3.5,-.3){$_{~4}$} \put(-2.8,-.3){$_{4}$}

\put(-1.8,0){$_{M(_s\gz)=}$}
\put(-.5,0.3){$_{~1}$}
\put(.1,0.3){$_{~1}$}
\put(-.5,0){$_{2~3~2~3~,}$}
\put(-.2,-.3){$_{~4}$} \put(.5,-.3){$_{4}$}
\put(-.75,-.3){$_{~4}$}

\put(1.7,0){$_{M(\gz_e)=}$}
\put(3.3,0.3){$_{~1}$}
\put(3,0){$_{3~2~3~2~,}$}
\put(3,-.3){$_{~4}$} \put(3.7,-.3){$_{4}$}

\put(4.9,0){$_{M(_s\gz_e)=}$}
\put(6.2,0.3){$_{~1}$}
\put(6.8,0.3){$_{~1}$}
\put(6.2,0){$_{2~3~2~3~2~.}$}
\put(6.5,-.3){$_{~4}$} \put(7.2,-.3){$_{4}$}
\put(5.95,-.3){$_{~4}$}
\end{pspicture}
\end{center}

Hence, we have an AR-sequence:
$$0\lra M(\gz)\lra M(_s\gz)\oplus M(\gz_e)\lra  M(_s\gz_e)\lra 0.$$
\end{example}
\medskip
\begin{rem}\label{rmk2}If $M(\gz)$ is an injective $J(Q,W)-$module, then it follows from the description in Remark \ref{rmk1}(2) that $_s\gz_e$ is an internal arc in $\ggz$ which implies $M(_s\gz_e)=0$. Dually, if $\gz \in \ggz$ is an internal arc, then $M(_s\gz_e)$ is a projective module in $\mod J(Q,W).$
\end{rem}

\subsection{Irreducible morphisms in $\mc_{(S,M)}$} Recall that $\mc:=\mc_\ggz$ is a triangulated Hom-finite $k$-category which is  2-Calabi-Yau.
We know that $\mc$ has a cluster-tilting object $T_\ggz=\tau_1\oplus\tau_2\oplus\cdots\oplus\tau_n$, where $\tau_1,\ldots\tau_n$ are the internal arcs in $\ggz$. Denote by $[1]$ the suspension functor of the triangulated category $\mc$. Moreover, since $
\mc$ is 2-Calabi-Yau, $\mc$ has AR-triangles and the AR-translation is given by  $\tau=[1]$ (see \cite{RV}).

Curves in $(S,M)$ which are non-contractible and not homotopic to a boundary arc give a parametrization of the isoclasses of string objects in $\mc$, we refer to them as non-boundary curves. We identify contractible curves and boundary arcs with the zero object in $\mc.$
We further say that a curve $\gz$ is not in $\ggz$ and write $\gz \not\in\ggz$ if it is non-contractible and not homotopic to an arc in the triangulation $\ggz$. Since the non-boundary arcs in $\ggz$ yield the cluster-tilting object $T_\ggz$,  the string module $M(\gz)$ is non-zero if $\gz \not\in\ggz$.
The following lemma from \cite{KR,KZ}  will be used frequently:

\begin{lem}\label{lem3}Let $\gz$ be a curve  in $(S,M)$ such that $\gz\not\in\ggz$. Then
 \begin{itemize}
 \item[(1)]$M(\gz)$ is  projective in $\mod J(Q,W)$ if and only if $\gz[1]\in \ggz,$
 \item[(2)]$M(\gz)$ is  injective in $\mod J(Q,W)$ if and only if $\gz[-1]\in \ggz.$
 \end{itemize}
 \end{lem}
From \cite{KZ}, we know that any sink (or source) map in $\mc$ is again a sink (or source) map in $\mod J(Q,W)$, therefore each AR-triangle in $\mc$ $$\gz'\lra\oplus\gz_i\lra\gz''\lra \gz'[1]$$
with $\gz'\not\in\ggz$ and $\gz'[-1]\not\in\ggz$ yields an AR-sequence in $\mod J(Q,W)$: $$0\lra M(\gz')\lra \oplus M(\gz_i)\lra M(\gz'')\lra 0.$$ Moreover, all AR-sequences in $\mod J(Q,W)$ are obtained in this way and the AR-translation in $\mod J(Q,W)$ is induced by that in $\mc$.
Combining this with Corollary \ref{ar-translation} and Remark \ref{rmk2}, we can easily get

\begin{prop}The AR-translation in $\mc_{(S,M)}$ is given by simultaneous counterclockwise rotation of starting and ending point of a curve to the next marked points at the boundary. That is
$\gz[-1]={_s\gz_e}$ if $\gz$ is a string object and
$(\lambda,b^n)[-1]=(\lambda,{_sb^n_e})$ if $(\lambda,b^n)$ is a band object in $\mc_{(S,M)}.$
\end{prop}
\begin{pf}Assume $\gz$ is a string object. If $\gz$ is a non-boundary curve in $(S,M)$ such that $M(\gz)$ is not injective, then
$M(\gz[-1])=M(\t^{-1}(\gz))=\t^{-1}(M(\gz))=M(_s\gz_e)$
which implies $\gz[-1]={_s\gz_e}.$ If $M(\gz)$ is injective, then Remark \ref{rmk2} shows that $_s\gz_e\in\ggz$ is an internal arc, and by Lemma \ref{lem3} and the discussion in Section 2 we obtain $\gz[-1]={_s\gz_e}.$
Similarly, if $\gz\in\ggz$ is an internal arc, then
$\gz[-1]={_s\gz_e}.$

Assume $(\lambda,b^n)$ is a band object where $\lambda\in k^*, b^n\in \Pi_1^*(S,M)$ with $n\geq 1$ and $b$ a closed curve, then
$M((\lambda,b^n)[-1])=M(\t^{-1}(\lambda,b^n))=\t^{-1}(M((\lambda,b^n)))\\=\t^{-1}(M(w(b),n,\phi_\lambda))=M(w(b),n,\phi_\lambda)
=M(w(_sb_e),n,\phi_\lambda)=M((\lambda, {_sb^n_e})),$ where $\phi_\lambda\in \mbox{Aut}(k^n)$ is induced by $\lambda.$ Hence $(\lambda,b^n)[-1]=(\lambda,{_sb^n_e})$ as a band object in $\mc_{(S,M)}.$
\end{pf}

We need the following lemma related to projective-injective modules in $\mod J(Q,W)$:

\begin{lem}\label{lem4}Let $\gz\not\in \ggz$ be a curve in $(S,M)$ such that $M(\gz_e)$ is a non-zero projective-injective module in $\mod J(Q,W)$.
Then there is a source map in $\mod J(Q,W)$:
$M({_s\gz_e})\lra M({_{s^2}\gz_e}).$
\end{lem}
\begin{pf}Assume $M(\gz_e)$ is a projective-injective module, then $w({\gz_e})$ must be an inverse (or a direct) string which both starts and ends on a peak and also both starts and ends in a deep. Without loss of generality, assume $w({\gz_e})=\bz_1^{-1}\cdots\bz_r^{-1}\bz_{r+1}^{-1}$ where $r\geq 0$.

\begin{center}
\begin{pspicture}(-3,-1.2)(3,1.2)

\put(-3,0){\circle*{.1}}
\put(3,0){\circle*{.1}}\put(-2,1){\circle*{.1}}\put(1.3,1){\circle*{.1}}
\put(2,-1){\circle*{.1}}\put(-.5,-1){\circle*{.1}}
\put(.5,-1){\circle*{.1}}\put(-2,-1){\circle*{.1}}
\psline[linewidth=1pt]{-}(2,-1)(3,0)(1.3,1)(-2,1)(-3,0)
\psline[linewidth=.5pt]{-}(3,0)(1.3,1)(-2,1)(-3,0)(1.3,1)(-2,-1)(-3,0)
\psline[linewidth=.5pt]{-}(-2,-1)(-.5,-1)(1.3,1)(.5,-1)(2,-1)(1.3,1)(3,0)(2,-1)
\psline[linewidth=1pt]{->}(2.8,.4)(2,.8)
\psline[linewidth=1pt]{<-}(-3,.4)(-2.5,.8)

\psline[linewidth=1pt,linestyle=dotted]{-}(.5,-1)(-.5,-1)
\psline[linewidth=.7pt,linecolor=blue]{->}(1.65,0)(.9,0)
\psline[linewidth=1pt,linestyle=dotted,linecolor=blue]{-}(.9,0)(0.4,0)
\psline[linewidth=.7pt,linecolor=blue]{->}(0.4,0)(-.35,0)
\psline[linewidth=.7pt,linecolor=blue]{->}(-.35,0)(-1.3,0.4)
\psline[linewidth=.5pt,linecolor=red]{-}(2,-1)(-2,1)
\psline[linewidth=.5pt,linecolor=red]{-}(3,0)(-2,1)
{\blue\uput[d](1.35,0){$_{\bz_1}$}
\uput[d](-.2,0.1){$_{\bz_r}$}
\put(-1.3,0){$_{\bz_{r+1}}$}}

\uput[r](2,-1){$_{b=e{(\gz)}}$}\uput[r](3,0){$_{c=e(\gz_e)}$}\uput[u](1.3,.9){$_{d}$}\uput[u](-2,.9){$_{a=s{(\gz)}}$}
\uput[l](-3,0){$_e$}\uput[l](-2,-1){$_f$}
\uput[d](0,-1){$_{Figure~6}$}
{\red\uput[d](1,-0.5){$_\gz$}\uput[d](2.2,.2){$_{\gz_e}$}}

\end{pspicture}
\end{center}

Let $a=s(\gz)=s(\gz_e), b=e(\gz), c=e(\gz_e)$ and $d,e,f$ be marked points as in Figure 6. Since $w({\gz_e})$ both starts and ends in a deep, $bc$ and $da$ are boundary arcs. Similarly, $ae$ and $dc$ are boundary arcs since $\gz_e$ both starts and ends on a peak, therefore $b,c,d,a,e$ lie consecutively on the same boundary.
\begin{itemize}
\item[(1)]If $M(_s\gz_e)$ is an injective module, that is $w({_s\gz_e})=\bz_1^{-1}\cdots\bz_r^{-1}$ both starts and ends on a peak, then $ef$ is a boundary arc. Then the definition of Pivot elementary move implies that $$M({{_{s^2}\gz}_{e}})=M(\bz_1^{-1}\cdots\bz_{r-1}^{-1})=M(_s\gz_e)/Soc(M(_s\gz_e))$$ and $_s\gz_{e^2}=de\in\ggz$ is an internal arc. Hence by Remark \ref{rmk1}(2) there is a source map: $M({_s\gz_e})\lra M({_{s^2}\gz_e}).$

\item[(2)]If $M(_s\gz_e)$ is not an injective module, that is $w({_s\gz_e})=\bz_1^{-1}\cdots\bz_r^{-1}$ does not start on a peak, then there exists $\az\in Q_1$ such that $w({_s\gz_e})\az$ is a string. By Theorem \ref{thm1}, there is an AR-sequence:
$$0 \lra M({_s\gz_e})\lra M({_{s^2}\gz_e})\oplus M(_s\gz_{e^2})\lra M({_{s^2}\gz_{e^2}})\lra 0,$$
where $w({_s\gz_e})=\bz_1^{-1}\cdots\bz_r^{-1}=U_\az$, $w({{_{s^2}\gz}_e})=w({_s\gz_e})_h=N_\az$ by Lemma \ref{lem1}, and $_s\gz_{e^2}=de\in\ggz$ is an internal arc in $\ggz$. Thus there is a source map: $M({_s\gz_e})\lra M({_{s^2}\gz_e})$ in $\mod J(Q,W)$.
\end{itemize}

\end{pf}

\begin{prop}\label{prop1}Let $\gz$ be a non-boundary curve and
$$0\lra M(\gz)\lra M(_s\gz)\oplus M(\gz_e)\lra M(_s\gz_e)\lra 0\qquad(3.1)$$
be an AR-sequence in $\mod J(Q,W)$ with $M(_s\gz)\neq 0\neq M(\gz_e)$, then there is an AR-triangle with two middle terms in $\mc$ as follows:
$$\gz\lra {_s\gz}\oplus\gz_e\lra {_s\gz_e}\lra \gz[1].$$
\end{prop}
\begin{pf}Since AR-sequences in $\mod J(Q,W)$ are induced by AR-triangles in $\mc$, we assume that
$$\gz\lra {_s\gz}\oplus\gz_e\oplus\sigma\lra {_s\gz_e}\lra \gz[1]\qquad(3.2)$$
is an AR-triangle in $\mc$ with $M(\sigma)=0$. It suffices to prove that $\sigma=0.$

Otherwise, $\sigma$ contains a direct summand $\rho$ which is an internal arc in $\ggz$. Then $_s\gz[1]\in \ggz$ or $\gz_e[1]\in \ggz$ since
$M(\rho[1])\neq 0$ in $\mod J(Q,W)$ by Lemma \ref{lem3} and since the number of summands in the middle term of an AR-sequence in $J(Q,W)$ is at most two.
Similarly, we conclude $_s\gz[-1]\in \ggz$ or $\gz_e[-1]\in \ggz$.
Hence there are 4 cases, we consider only two of them, the others are similar.

\begin{itemize}
\item[(1)]Assume $_s\gz[1]\in \ggz$ and $\gz_e[-1]\in \ggz$, then Lemma \ref{lem3} implies that $M(_s\gz)$ is projective and $M(\gz_e)$ is injective. Hence $M(\gz)$ must be a direct summand of rad$M(_s\gz)$ and $M(_s\gz_e)$ must be a direct summand of $M(\gz_e)/$Soc$M(\gz_e)$.
This yields a contradiction by comparing the dimensions in $(3.1)$.
\item[(2)]Assume $\gz_e[1]\in \ggz$ and $\gz_e[-1]\in \ggz$, then Lemma \ref{lem3} implies that $M(\gz_e)$ is a projective-injective module, and so by Lemma \ref{lem4} there is a source map
    $$M({_s\gz_e})\lra M({_{s^2}\gz_e})\qquad(3.3)$$
    in $\mod J(Q,W)$. But after shifting $(3.2)$ by [-1], there is an AR-sequence in $\mod J(Q,W)$:
     $$0\lra M(_s\gz_e)\lra M({_{s^2}\gz_e})\oplus M(_s\rho_e)\lra M(_{s^2}\gz_{e^2})\lra 0.$$
     with $M(_s\rho_e)=M(\rho[-1])\neq 0$ a projective module. This contradicts the fact that the source map $(3.3)$ ends with just one indecomposable module.  \end{itemize}
\end{pf}

Since neither injective modules nor projective modules occur in homogeneous tubes, the following corollary related to band objects can be obtained similarly.
\begin{cor}\label{bands}
Let $(\lambda,b^n)$ be a band object in $\mc$ where $\lambda\in k^*, b^n\in \Pi_1^*(S,M)$ with $n \geq 1$ and $b$ itself is not a proper power of an element in $\Pi_1^*(S,M)$, then there is an AR-triangle in $\mc:$
$$(\lambda,b^n)\lra(\lambda,b^{n+1})\oplus (\lambda, b^{n-1})\lra (\lambda, b^n)\lra(\lambda, b^n)[1],$$
where $(\lambda, b^0)$ is the zero object in $\mc$.
\end{cor}
Therefore, the band objects are closed under irreducible morphisms and the corresponding AR-components in $\mc$ are homogeneous tubes. It remains to consider the string objects in $\mc.$

\begin{lem}\label{lem6}Let $\gz$ be an internal arc in $\ggz$, then the AR-triangle in $\mc$ starting in $\gz$ is of the form
$$\gz\lra {_s\gz}\oplus\gz_e\lra{_s\gz_e}\lra\gz[1].$$
\end{lem}
\begin{pf}
Suppose $\gz$ is an internal arc in $\ggz$, then $M(_s\gz_e)$ is projective in $\mod J(Q,W)$ by Lemma \ref{lem3}. We consider $\gz[-1]={_s\gz_e}$.
\begin{itemize}
\item[(1)]If $M(_s\gz_e)$ is not injective, then Theorem \ref{thm1} induces an AR-sequence in $\mod J(Q,W)$:
$$\qquad\qquad 0\lra M(_s\gz_e)\lra M({_{s^2}\gz_e})\oplus M(_s\gz_{e^2})\lra M(_{s^2}\gz_{e^2})\lra 0\qquad(3.4)$$
The definition of a cluster-tilting object guarantees that  $_{s^2}\gz_e$ and $_s\gz_{e^2}$ cannot be internal arcs in $\ggz$. If neither $_{s^2}\gz_e$ nor $_s\gz_{e^2}$ is boundary arc (thus $M(_{s^2}\gz_e)\neq 0 \neq M(_s\gz_{e^2})$), then Proposition \ref{prop1} induces an AR-triangle
$$_s\gz_e \lra {_{s^2}\gz_e}\oplus {_s\gz}_{e^2} \lra {_{s^2}\gz}_{e^2}\lra {_s\gz_e}[1]$$
which yields the AR-triangle
$$\gz\lra {_s\gz}\oplus\gz_e\lra {_s\gz_e}\lra \gz[1].$$
 Assume now that one of $_{s^2}\gz_e$ and $_s\gz_{e^2}$ is a boundary arc.
 Without loss of generality, let $_{s^2}\gz_e$ be a boundary arc, then by definition of a cluster-tilting object, $\Hom_\mc(_s\gz_e,\sigma)=0$ for any $\sigma\in\ggz$, hence $(3.4)$ induces an AR-triangle of the form
 $$_s\gz_e \lra {_{s^2}\gz_e}\oplus {_s\gz}_{e^2} \lra {_{s^2}\gz}_{e^2}\lra {_s\gz_e}[1]$$
 with ${_{s^2}\gz_e}=0$. The definition of a Pivot elementary move implies that $\gz_e$ is also a boundary arc, hence there is an AR-triangle of the form
$$\gz\lra {_s\gz}\oplus\gz_e\lra {_s\gz_e}\lra \gz[1]$$ where $\gz_e=0$ is a boundary arc in $\mc$.
\medskip

\item[(2)]If $M(_s\gz_e)=M(\gz[-1])$ is injective, then $M(_s\gz_e)$ is a projective-injective module by Lemma \ref{lem3}.
By definition of a cluster-tilting object and Remark \ref{rmk1}(3) we can assume that the AR-triangle starting in $\gz$ is of the form
$$\gz\lra \da_1\oplus \da_2\lra {_s\gz_e}\lra\gz[1]$$
where neither $\da_1$ nor $\da_2$ are internal arcs in $\ggz$.
Since $M(_s\gz_e)$ is a projective-injective module in $\mod J(Q,W)$, one of  $\da_1$ and $\da_2$ must be a boundary arc. Without loss of generality, assume $\da_2$ is a boundary arc, then Lemma \ref{lem4} implies $M(\da_1)=\rad M(_s\gz_e)=M({_s}\gz)$, that is $\da_1={_s}\gz$ and $\da_2=\gz_e$ is a boundary arc. Therefore the AR-triangle is of the form
$$\gz\lra {_s\gz}\oplus\gz_e\lra {_s\gz_e}\lra \gz[1]$$ where $\gz_e=0$ is a boundary arc.
\end{itemize}
\end{pf}
As shown in Remark \ref{rmk1}(3), the middle terms of AR-sequences in $\mod J(Q,W)$ contain at most two indecomposable summands. The following proposition establishes the same result for AR-triangles in $\mc$.
\begin{prop}\label{prop2}The number of indecomposable summands in the middle term of an AR-triangle in $\mc_{(S,M)}$ is  at most two.
\end{prop}
\begin{pf}
We only need to consider string objects in $\mc$. Let therefore $\gz$ be a non-boundary curve, and
$$\gz\lra \bigoplus_{i=1}^r \da_i\oplus\bigoplus_{j=1}^s \t_j\lra\gz[-1]\lra\gz[1]\qquad(3.5)$$
be an AR-triangle in $\mc$ starting in $\gz$ with $\da_i\not\in\ggz$ and where $\t_j$ are internal arcs in $\ggz$.
Remark \ref{rmk1}(3) implies $r\leq 2$ and $s\leq 2.$ If $\gz\in \ggz$ or $\gz[-1]\in\ggz$, then by Lemma \ref{lem6}, there is nothing to prove. We suppose now $\gz\not
\in \ggz$ and $\gz[-1]\not\in\ggz.$ If $r=2$, Proposition \ref{prop1} implies $s=0$. If $r=1$ and $s=2$, then $(3.5)$ induces an AR-sequence with one middle term in $\mod J(Q,W)$:
$$0\lra M(\gz)\lra M(\da_1)\lra M(\gz[-1])\lra 0,$$ hence there exists $\az\in Q_1$ such that $w({\da_1})=N_\az$ by Remark \ref{rmk1}(1). On the other hand, if we shift $(3.5)$ by [1] and [-1], then Remark \ref{rmk1}(3) implies that $\da_1[1]$ and $\da_1[-1]$ must be the internal arcs in $\ggz$ since the $M(\t_1[1])\neq 0\neq M(\t_2[1])$ and $M(\t_1[-1])\neq 0\neq  M(\t_2[-1])$ in $\mod J(Q,W)$ by Lemma \ref{lem3}.
Thus $M(\da_1)=M(N_\az)$ is a projective-injective module in $\mod J(Q,W)$. The definition of $N_\az$ implies that the quiver $Q$ is of type $A_2.$ But it is impossible since the cluster category of type $A_2$ is well-known, and there are no AR-triangle with more than two summands in the middle term.
\end{pf}
Now we consider the AR-sequences with just one middle term in $\mod J(Q,W)$.

\begin{lem}\label{lem5} Let $\gz$ be a non-boundary curve in $(S,M)$ such that $w(\gz)=U_\az,$ for some $\az\in Q$. Then $w({_s\gz})=N_\az$ and $w({\gz_e})=0$ is a zero string, and the AR-triangle in $\mc$ starting in $\gz$ is of the form:
$$\gz\lra {_s\gz}\oplus\gz_e\lra {_s\gz_e}\lra \gz[1].$$

\end{lem}
\begin{pf}Since $w(\gz)=U_\az$ does not start but ends on a peak, we know by Lemma \ref{lem1} that $w({_s\gz})={w(\gz)}_h=N_\az$ and by Lemma \ref{lem2} that $w({\gz_e})={_cw}(\gz)={_cU_\az}=0$ is a zero string. Hence $\gz_e$  might be an internal arc in $\ggz$ or a boundary arc in $(S,M)$, see Figure 7

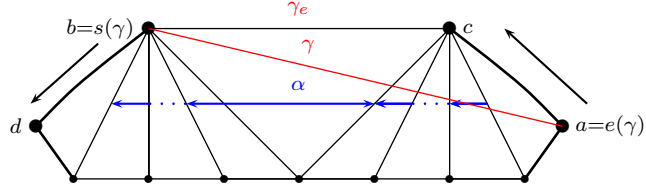
\begin{center}
\begin{pspicture}(-3,-1.2)(3,1.3)
\put(-3.5,-.3){\circle*{.17}}
\put(3.5,-.3){\circle*{.17}}\put(-2,1){\circle*{.17}}\put(2,1){\circle*{.17}}
\put(0,-1){\circle*{.1}}
\put(-3,-1){\circle*{.1}}\put(-2,-1){\circle*{.1}}\put(-1,-1){\circle*{.1}}
\put(3,-1){\circle*{.1}}\put(2,-1){\circle*{.1}}\put(1,-1){\circle*{.1}}
\pscurve[linewidth=1pt]{-}(3.5,-.3)(3,.2)(2,1)
\pscurve[linewidth=1pt]{-}(-3.5,-.3)(-3,.2)(-2,1)
\psline[linewidth=1pt]{-}(3.5,-.3)(3,-1)
\psline[linewidth=1pt]{-}(-3.5,-.3)(-3,-1)
\psline[linewidth=.5pt]{-}(-3.5,-.3)(-3,-1)(3,-1)(3.5,-.3)
\psline[linewidth=.5pt]{-}(-3,-1)(-2,1)(-2,-1)(-2,1)(-1,-1)(0,-1)(-2,1)(2,1)(0,-1)(1,-1)(2,1)(2,-1)(3,-1)(2,1)
\psline[linewidth=.8pt,linecolor=blue]{->}(2.5,0)(2,0)
\psline[linewidth=.8pt,linecolor=blue]{->}(1.5,0)(1,0)
\psline[linewidth=1pt,linecolor=blue,linestyle=dotted]{-}(2,0)(1.5,0)
\psline[linewidth=.8pt,linecolor=blue]{<-}(-2.5,0)(-2,0)
\psline[linewidth=.8pt,linecolor=blue]{<-}(-1.5,0)(-1,0)
\psline[linewidth=1pt,linecolor=blue,linestyle=dotted]{-}(-2,0)(-1.5,0)
\psline[linewidth=.8pt,linecolor=blue]{->}(2.5,0)(2,0)
\psline[linewidth=.8pt,linecolor=blue]{->}(1.5,0)(1,0)
\psline[linewidth=.8pt,linecolor=blue]{->}(-1,0)(1,0)
{\blue\uput[u](0,0){$_{\az}$}}
\psline[linewidth=.5pt,linecolor=red]{-}(3.5,-.3)(-2,1)
\uput[r](3.5,-.3){$_{a=e{(\gz)}}$}\uput[r](2,1){$_c$}
\uput[l](-2,1){$_{b=s{(\gz)}}$}\uput[l](-3.5,-.3){$_{d}$}
\uput[d](0,-1){$_{Figure~7}$}{\red\uput[u](0,1){$_{\gz_e}$}}
{\red\uput[u](0,.5){$_\gz$}}
\psline{->}(3.7,0)(2.6,1)\psline{<-}(-3.7,0)(-2.8,.8)
\end{pspicture}
\end{center}

Let $b=s(\gz), a=e(\gz)$ be the endpoints of $\gz$ and let $d$ be the marked point lying clockwise next to $b$ on the same boundary. Since $U_\az$ ends on a peak, $ac$ is a boundary arc. The definition of Pivot elementary move implies $c=e(\gz_e)$, hence there is an AR-sequence with one middle term in $\mod J(Q,W)$ (See Remark \ref{rmk1}(1)):
$$0\lra M(\gz)\lra M(_s\gz)\lra M(_s\gz_e)\lra 0.\qquad(3.6)$$
\begin{itemize}
\item[(1)]If $\gz_e\in\ggz$, then Lemma \ref{lem3} implies that $M(\gz_e[1])$ is an injective module. Remark \ref{rmk1}(2) induces an irreducible morphism $M(\gz_e[1])\lra M(\gz)$ in $\mod J(Q,W)$, hence there are two irreducible morphisms in $\mc$: $\gz_e[1]\lra \gz$ and $\gz\lra \gz_e[1][-1]=\gz_e$. Combining this with $(3.6)$ and Proposition \ref{prop2}, we get the AR-triangle in $\mc:$
   $$\gz\lra {_s\gz}\oplus\gz_e\lra {_s\gz_e}\lra \gz[1].$$
\item[(2)]If  $\gz_e$ is a boundary arc, then $w(\gz)$ starts in a deep and ends on a peak which implies $M(\gz)=M(U_\az)$ is projective in $\mod J(Q,W)$, thus $\gz[1]\in \ggz$ by Lemma \ref{lem3}. By $(3.6)$ and the definition of a cluster-tilting object, the AR-triangle starting in $\gz$ is of the form
    $$\gz\lra {_s\gz}\oplus\gz_e\lra {_s\gz_e}\lra \gz[1]$$
    where $\gz_e=0$ in $\mc.$
\end{itemize}
\end{pf}

By the following theorem we finally obtain  the results formulated in Theorem \ref{irred} in the introduction:
\begin{thm}\label{thm2}Let $\ggz$ be a triangulation of a marked surface $(S,M)$, and let $\gz$ be a non-boundary curve in $(S,M)$.
Then each irreducible morphism in $\mc_{(S,M)}$ starting in $\gz$ is obtained by Pivot elementary moves on endpoints of $\gz$. The AR-triangle starting in $\gz$ is of the form:
$$\gz\lra {_s\gz}\oplus\gz_e\lra {_s\gz_e}\lra \gz[1].$$
\end{thm}
\begin{pf} Lemma \ref{lem6} implies the case when $\gz$ or $\gz[-1]={_s\gz_e}$ is an internal arc in $\ggz$. Proposition \ref{prop1} and Lemma \ref{lem5} yield the remaining case.
\end{pf}

\subsection{AR-components in $\mc_{(S,M)}$}
The aim of this subsection is to describe the AR-components in  $\mc_{(S,M)}$.
From Corollary \ref{bands} we know that all band objects of $\mc_{(S,M)}$ lie in homogeneous tubes, so we focus on string objects from now on.
It follows from Theorem \ref{thm2} that each string object $\gz$ in $\mc_{(S,M)}$ is starting point of a mesh with two middle terms in the AR-quiver of   $\mc_{(S,M)}$ except when one of ${_s\gz}$ or $\gz_e$ is a boundary arc.
The situation where one of ${_s\gz}$ or $\gz_e$ is a boundary arc is explicitly described in the following corollary.

\begin{cor}\label{boundary}
Let $C$ be a boundary component of $S$ with $t$ marked points. If we choose a numeration and an orientation of the boundary arcs $\delta^1, \ldots, \delta^t$ such that all $\delta^i$ are clockwise oriented and $e(\delta^i) = s(\delta^{i+1})$ for all $i$, then the objects $\delta^i_{e^j}$ with $i=1,\ldots,t$ and  $j \ge 1$ form a tube of rank $t$ in $\mc_{(S,M)}$.
Moreover, there is a bijection between the boundary components of $S$ and the tubes in  $\mc_{(S,M)}$ which are not formed by band objects.
\end{cor}
\begin{pf} The proof follows easily from the description of the AR-triangles in Theorem \ref{thm2} once the orientation of the boundary arcs is chosen, we illustrate the situation where $t=3$ in following figure.
Note that $\delta^i_{e^{t-1}}$ is a non-contractible closed curve for all surfaces except a disc, where all closed curves are contractible. Thus, for a disc the objects $\delta^i_{e^{t-1}}$ are zero, and the remaining objects $\delta^i_{e^j}$ with  $j > t-1$ have to be identified accordingly, see the description in \cite{CCS}.
\end{pf}

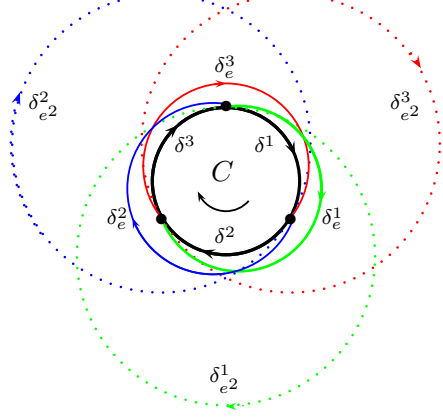
\begin{center}
\begin{pspicture}(-4,-3)(4,2.2)
\put(-.2,0){$C$}
\pscircle[linewidth=1.2pt](0,0){1}
\psarc[linewidth=.7pt]{<-}(0,0){.4}{-160}{-41}
\psarc[linewidth=1.2pt]{<-}(0,0){.98}{-110}{-41}
\psarc[linewidth=1.2pt]{<-}(0,0){.98}{-230}{-110}
\psarc[linewidth=1.2pt]{<-}(0,0){.98}{-340}{-230}
\uput[u](0,-1){$_{\da^2}$}\uput[r](-.86,.5){$_{\da^3}$}\uput[l](.86,.5){$_{\da^1}$}
\uput[u](0,-3){$_{\da^1_{e^2}}$}\uput[u](0,1.2){$_{\da^3_{e}}$}
\uput[l](2.8,1){$_{\da^3_{e^2}}$}\uput[r](-2.8,1){$_{\da^2_{e^2}}$}
\uput[r](1.1,-.5){$_{\da^1_e}$}\uput[l](-1.1,-.5){$_{\da^2_e}$}
\pscircle[linewidth=1pt,linestyle=dotted,linecolor=green](0,-1){2}
\psarc[linewidth=.7pt,linestyle=dotted,linecolor=green]{<-}(0,-1){2}{-90}{-80}
\pscircle[linewidth=1pt,linecolor=red,linestyle=dotted](0.86,.5){2}
\psarc[linewidth=1pt,linecolor=red,linestyle=dotted]{<-}(0.86,.5){2}{30}{40}
\pscircle[linewidth=1pt,linecolor=blue,linestyle=dotted](-0.86,.5){2}
\psarc[linewidth=1pt,linecolor=blue,linestyle=dotted]{<-}(-0.86,.5){2}{160}{-140}
\psarc[linewidth=.7pt,linecolor=red](0,0.2){1.1}{-41}{-138}
\psarc[linewidth=.7pt,linecolor=red]{<-}(0,.2){1.1}{90}{98}
\psarc[linewidth=1pt,linecolor=green]{<-}(.17,-0.1){1.1}{-10}{61}
\psarc[linewidth=1pt,linecolor=green](.17,-0.1){1.1}{-161}{98}
\psarc[linewidth=.7pt,linecolor=blue](-.172,-0.1){1.14}{82}{-23}
\psarc[linewidth=.7pt,linecolor=blue]{<-}(-.172,-0.1){1.14}{200}{-82}
\put(0,1){\circle*{.15}}
\put(-.86,-.5){\circle*{.15}}
\put(.86,-.5){\circle*{.15}}
\end{pspicture}
\end{center}

The preceding corollary yields a description of all AR-components (formed by string objects) which contain some meshes with only one middle term. The remaining AR-components in  $\mc_{(S,M)}$ are formed by meshes with exactly two middle terms. If $\Lambda$ is one such component, one can choose one object $\gz$ in $\Lambda$ and then the component is formed by all $_{s^i}\gz_{e^j}$ with $i,j \ge 1$.
In case there are no identifications of the objects, one obtains thus components of the form  $\mathbb{Z}\mathbb{A}^\infty_\infty$.
If $S$ is neither a disc nor an annulus, then the Jacobian algebra associated to a triangulation will in general be of non-polynomial growth (see \cite{ABCP}), and there are plenty of components of the form  $\mathbb{Z}\mathbb{A}^\infty_\infty$ in the AR-quiver of $\mc_{(S,M)}$.
If $S$ is an annulus, the cluster category $\mc_{(S,M)}$ is of type $\widetilde{\mathbb{A}}$ and there is no component of the form $\mathbb{Z}\mathbb{A}^\infty_\infty$.

\section{Flips and Mutations}

In the previous sections we fixed one triangulation $\Gamma$ of $(S,M)$ and studied the irreducible morphisms and AR-components of the cluster category $\mc_\ggz$ defined by  the quiver with potential $(Q_\ggz,W_\ggz)$. The aim of this section is to study the effect of a change in the triangulation on the cluster category $\mc_\ggz$.

If $\tau_i$ is an internal arc in $\ggz$, then there exists exactly one internal arc $\tau_i'\neq\tau_i$ in $(S,M)$ such that $f_{\tau_i}(\ggz):=(\ggz\backslash\{\tau_i\})\cup\{\tau_i'\}$ is also a triangulation of $(S,M)$.
In fact, the internal arc ${\tau_i}$ is a diagonal in the quadrilateral formed by the two triangles of $\ggz$ containing ${\tau_i}$, and ${\tau'_i}$ is the other diagonal in that quadrilateral, see \cite{FST}.
We denote $\tau_i'$ by $f_\ggz(\tau_i)$ and say that $f_{\tau_i}(\ggz)$ is obtained from $\ggz$ by applying a flip along $\tau_i$.
In fact, by applying iterated flips one can obtain all triangulations of $(S,M)$:

\begin{thm}[\cite{FST}]\label{fst} For any two triangulations of $(S,M)$ there is a sequence of flips which transforms one triangulation into the other.
\end{thm}

As shown in Theorem 1.1, we can view the non-boundary curves in $(S,M)$  as objects in $\mc_\ggz$.
If we denote  all internal arcs of $\ggz$ by $\tau_1, \ldots,\tau_n$, then their direct sum $T_{\ggz}=\tau_1\oplus\tau_2\cdots\oplus\tau_n$ is a cluster-tilting object in $\mc_\ggz$. The following theorem is adapted from \cite{IY} to our setup.

\begin{thm}[\cite{IY}] Let $\tau_i\in\ggz$ be an internal arc, then  there is a curve $\tau^*_i$ in $(S,M)$ (unique up to homotopy) which is not homotopic to $\tau_i$ such that the object $\mu_{\tau_i}(T_\ggz)$ obtained from $T_\ggz$ by replacing $\tau_i$ with $\tau^*_i$ is also a cluster-tilting object in $\mc_\ggz$.
\end{thm}

The object  $\mu_{\tau_i}(T_\ggz)$ is called the \emph{mutation} of $T_\ggz$ in $\tau_i$, and $(\tau_i,\tau^*_i)$ is called an \emph{exchange pair} in $\mc_\ggz$.
As shown in \cite{IY}, any exchange pair $(\tau_i, \tau^*_i)$ induces the following non-split triangles (unique up to isomorphism) which are referred to as  \emph{exchange triangles}:
$$\tau_i\overset{f}{\lra} \tau\lra\tau^*_i\lra\tau_i[1] \mbox{\quad and \quad } \tau_i^*\lra \tau\overset{g}{\lra}\tau_i\lra\tau_i^*[1]$$
Here $f$ is a minimal left add$(T_{\ggz}\setminus\{\tau_i\})$-approximation and $g$ is a minimal right add$(T_{\ggz}\setminus\{\tau_i\})$-approximation.
Since $$\End(T_\ggz,T_\ggz)\cong \End(T_\ggz[1],T_\ggz[1])\cong J(Q_\ggz,W_\ggz),$$ we know that the quiver of the endomorphism algebra of $T_\ggz$  in $\mc_\ggz$ does not contain loops at any vertex. Hence we obtain the following lemma which is a special case of Lemma 7.5 in \cite{K1}.

\begin{lem}\label{exch}Let $\tau_i$ be an internal arc of $\ggz$, then the exchange triangles are given by
$$\tau_i\overset{f}{\lra} \bigoplus\limits_{\tau_j\ra\tau_i}\tau_j\lra\tau^*_i\lra\tau_i[1] \mbox{\quad and \quad} \tau_i^*\lra \bigoplus\limits_{\tau_i\ra\tau_k}\tau_k \overset{g}{\lra}\tau_i\lra\tau_i^*[1]$$
where $f$ is a minimal left $\add (T_\ggz\setminus \tau_i)$-approximation and $g$ is a minimal right $\add(T_\ggz\setminus\tau_i)$-approximation.
\end{lem}

For an internal arc $\t_i\in\ggz$, we discussed the definition of  the flip $f_\ggz(\t_i)$ above. On the other hand, if we view $\t_i$ as an indecomposable rigid object in $\mc_\ggz$, also the mutation $\mu_\ggz(\tau_i):=\t_i^*$ of $\t_i$ is defined. The following theorem shows that flip and mutation of an internal arc are compatible (viewed as objects in the cluster category $\mc_\ggz$):

\begin{thm}\label{thm3} Let $\ggz$ be a triangulation of $(S,M)$ and let  $\tau_i\in\ggz$ be an internal arc.
Then
$$\mu_\ggz(\tau_i)=f_\ggz(\tau_i).$$
\end{thm}
\begin{pf}The internal arc $\tau_i$ is a diagonal of a  quadrilateral  formed by the internal arcs $\tau_{j_1},\tau_{j_2},\tau_{k_1},\tau_{k_2}$ as follows:

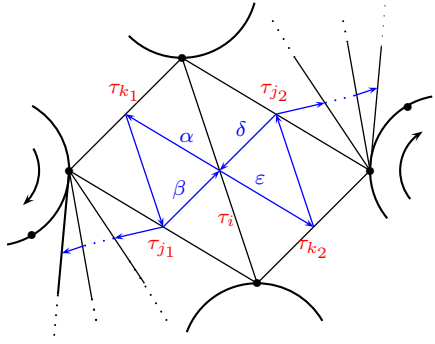
\begin{center}
\begin{pspicture}(-7,-2.2)(7,1.7)
\psarc[arcsepB=.5pt]{-}(-3,0){1}{-80}{-300}\psarc[arcsepB=.5pt]{<-}(-3,0){.6}{-50}{-330} \put(-2.5,-.85){\circle*{.1}}
\psarc[arcsepB=.5pt]{-}(-.5,2.5){1}{225}{-15}
\psarc[arcsepB=.5pt]{-}(.5,-2.5){1}{30}{-200}
\psline[linewidth=.5pt]{-}(-.5,1.5)(-2,0)(.5,-1.5)(2,0)(-.5,1.5)(.5,-1.5)
\put(-2,0){\circle*{.1}} \put(2,0){\circle*{.1}} \put(-.5,1.5){\circle*{.1}} \put(.5,-1.5){\circle*{.1}}
\psline[linewidth=.5pt,linecolor=blue]{->}(0,0)(-1.25,.75)
{\blue\uput[l](.55,.55){$_{\da}$}\uput[l](-.25,-.25){$_{\bz}$}\uput[l](-.15,.45){$_{\az}$}
\uput[u](.55,-.35){$_{\vz}$}}{\red\uput[d](1.25,-.75){$_{\tau_{k_2}}$}
\uput[d](-.75,-.75){$_{\tau_{j_1}}$}\uput[u](-1.25,.75){$_{\tau_{k_1}}$}\uput[u](.75,.75){$_{\tau_{j_2}}$}}
\psline[linewidth=.5pt,linecolor=blue]{<-}(1.25,-.75)(0,0)
\psline[linewidth=.5pt,linecolor=blue]{->}(.75,.75)(0,0)
\psline[linewidth=.5pt,linecolor=blue]{->}(-1.25,.75)(-.75,-.75)
\psline[linewidth=.5pt,linecolor=blue]{->}(-.75,-.75)(0,0)
\psline[linewidth=.5pt,linecolor=blue]{->}(1.25,-.75)(.75,.75)
\psline[linewidth=.5pt]{-}(-2,0)(-1,-1.5){\red\uput[d](0.1,-.4){$_{\tau_i}$}}
\psline[linewidth=.5pt]{-}(-2,0)(-1.7,-1.7)
\psline[linewidth=.5pt]{-}(2,0)(1,1.5)\psline[linewidth=.5pt,linestyle=dotted]{-}(2.2,2.2)(2.15,1.65)\psline[linewidth=.5pt]{-}(2,0)(2.15,1.65)
\psline[linewidth=.5pt]{-}(2,0)(1.7,1.7)
\psline[linewidth=.8pt,linestyle=dotted]{-}(-.7,-2)(-1.2,-1.2)\psline[linewidth=.8pt,linestyle=dotted]{-}(.7,2)(1,1.5)
\psline[linewidth=.8pt,linestyle=dotted]{-}(-1.7,-1.7)(-1.64,-2)\psline[linewidth=.8pt,linestyle=dotted]{-}(1.7,1.7)(1.64,2)
\psline[linewidth=.8pt,linestyle=dotted]{-}(-1.7,-1.7)(-1.64,-2)
\psline[linewidth=.8pt,linestyle=dotted]{-}(-2.2,-2.2)(-2.15,-1.65)\psline[linewidth=.8pt]{-}(-2,0)(-2.15,-1.65)
\psarc[arcsepB=.5pt]{-}(3,0){1}{100}{220}\psarc[arcsepB=.5pt]{<-}(3,0){.6}{120}{220}\put(2.5,.85){\circle*{.1}}

\psline[linewidth=.5pt,linecolor=blue]{->}(-.75,-.75)(-1.4,-.9)\psline[linewidth=.8pt,linestyle=dotted,linecolor=blue]{-}(-1.4,-.9)(-1.8,-1)
\psline[linewidth=.5pt,linecolor=blue]{->}(-1.8,-1)(-2.1,-1.1)
\psline[linewidth=.5pt,linecolor=blue]{->}(.75,.75)(1.4,.9)\psline[linewidth=.8pt,linestyle=dotted,linecolor=blue]{-}(1.4,.9)(1.8,1)
\psline[linewidth=.5pt,linecolor=blue]{->}(1.8,1)(2.1,1.1)

\end{pspicture}
\end{center}

The triangles of $\ggz$ containing $\tau_i$ induce arrows $\az,\bz,\da,\varepsilon$ in $Q_\ggz$ as indicated above.
By Lemma \ref{exch} there is a non-split triangle in $\mc_\ggz$
$$\tau_i\overset{f}{\lra} \tau_{j_1}\oplus\tau_{j_2}\lra\tau^*_i\lra\tau_i[1]$$
where $f$ is a minimal left add$(T_\ggz\setminus\tau_i)$-approximation. We obtain the following right exact sequence in $\mod J(Q_\ggz,W_\ggz)$ by applying $M(-)\cong\Ext_{\mc_\ggz}^1(T_\ggz,-)$ to this triangle:
$$M(\tau_i[-1])\lra M(\tau_{j_1}[-1])\oplus M(\tau_{j_2}[-1])\lra M(\tau_i^*[-1])\lra 0$$
This sequence is in fact a minimal projective resolution in $\mod J(Q_\ggz,W_\ggz)$ whose projective modules can be described by strings as follows:
$$M(\tau_i[-1])=M({V_\az}^{-1}V_\varepsilon),$$
$$M(\tau_{j_1}[-1])=M({V_\az}^{-1}\bz V_\bz) \mbox{\quad and \quad} M(\tau_{j_2}[-1])=M({V_\da}^{-1}\da^{-1} V_\varepsilon).$$ Hence $M(\tau_i^*[-1])\cong M(V_\da^{-1}\da^{-1}\bz V_\bz)$, which implies $w(_s({\tau_i^*})_e)=V_\da^{-1}\da^{-1}\bz V_\bz$ and thus $\mu_\ggz{(\tau_i)}=\tau_i^*=f_\ggz(\tau_i).$
\end{pf}

For any curve $\gz$ in $(S,M)$, the definition of the string $w(\gz)$ depends on the Jacobian algebra $J(Q_\ggz,W_\ggz)$. In order to compare string modules in two Jacobian algebras arising from different triangulations of $(S,M)$,  we denote in the following by $w(\ggz,\gz)$ and $M(\ggz,\gz)$ the string and the string module in the Jacobian algebra $J(Q_\ggz,W_\ggz)$. Similarly, the band in $J(Q_\ggz,W_\ggz)$ given by a closed curve $b$ is denoted by $w(\ggz,b).$ If $\gz=\t_i\in\ggz,$ we denote by $M(\ggz,\t_i)$ the associated simple decorated representation of $(Q_\ggz,W_\ggz)$ (see more details in \cite{DWZ,LF2}). It is shown in \cite{LF2} that the flips of triangulations are compatible with the mutations of decorated representations.

On the other hand, let $\ggz'=f_{\t_i}(\ggz)$ be the triangulation of $(S,M)$ obtained by a flip along $\t_i$, then \cite{KY} establishes an equivalence $\bar{\mu}_i:\mc_\ggz\lra\mc_{\ggz'}$ for each $1\leq i\leq n.$ Viewing the indecomposable objects in $\mc_\ggz$ as curves and closed curves in $(S,M)$, or as indecomposable decorated representation of $(Q_\ggz,W_\ggz),$ then \cite{P} shows the compatibility between the equivalence $\bar{\mu}_i$ and the mutation of decorated representations of $(Q_\ggz,W_\ggz).$ Therefore, we get following lemma.

\begin{lem}\label{lem7}Let $\ggz'=f_{\t_i}(\ggz)$ be the triangulation of $(S,M)$ obtained by a flip along $\t_i$ and let $T_{\ggz'}$ be the corresponding cluster-tilting object in $\mc_{\ggz'}$. If
$\bar{\mu}_i:\mc_\ggz\lra\mc_{\ggz'}$ denotes the equivalence from \cite{KY}, then
$$\Ext^1_{\mc_{\ggz'}}(T_{\ggz'},\bar{\mu}_i(\gz))=M(\ggz',\gz)$$
for any non-boundary curve $\gz$ in $(S,M).$
\end{lem}
\begin{pf}
Let $\mbox{Ind}_s\mc_\ggz$ be the set of all indecomposable string objects in $\mc_\ggz$ which is indexed by non-boundary curves in $(S,M)$, and $$D_\ggz=\{M(\ggz,\gz)|\ \gz\ curve\ in\ (S,M)  \}.$$ Then by Proposition 4.1 in \cite{P}, we have following commutative diagram
\begin{center}
\begin{pspicture}(-1.7,-.7)(1.2,1)
\put(-2,.7){${\mbox{Ind}_s\mc_\ggz}$}
\put(-2,-.7){${\mbox{Ind}_s\mc_{\ggz'}}$}
\put(1,.7){${D_\ggz}$}
\put(1,-.7){${D_{\ggz'}}$}
\psline[linewidth=.7pt]{->}(-.8,0.7)(.8,.7)
\psline[linewidth=.7pt]{->}(-.8,-0.7)(.8,-.7)
\psline[linewidth=.7pt]{<-}(-1.7,-0.3)(-1.7,.5)
\psline[linewidth=.7pt]{<-}(1.2,-0.3)(1.2,.5)
{\tiny\uput[u](0,.7){$_{\Ext^1_{\mc_\ggz}(T_\ggz,-)}$}}
{\tiny\uput[u](0,-.7){$_{\Ext^1_{\mc_{\ggz'}}(T_{\ggz'},-)}$}}
\uput[r](1.2,0){$_{\mu_i}$}
\uput[l](-1.7,0){$_{\bar{\mu}_i}$}
\end{pspicture}
\end{center}
where $\mu_i$ is the mutation of decorated representations. Therefore, for each curve $\gz$ in $(S,M)$,
$$\Ext^1_{\mc_{\ggz'}}(T_{\ggz'},\bar{\mu}_i(\gz))=\mu_i(\Ext^1_{\mc_\ggz}(T_\ggz,\gz))$$
$$\cong \mu_i(M(\ggz,\gz))\cong M(f_{\t_i}(\ggz),\gz)=M(\ggz', \gz)$$
where the last isomorphism is given by the main result in \cite{LF2}. This completes the proof.
\end{pf}
Since $\bar{\mu}_i:\mc_{\ggz}\ra\mc_{\ggz'}$ is an equivalence, the objects $\bar{\mu}_i(\gz)$ in $\mc_{\ggz'}$ can also be described by curves in $(S,M)$ and we denote $\bar{\mu}_i(\gz)$ again by $\gz$ in $\mc_{\ggz'}$.

\begin{cor}\label{cor1}Each triangulation of $(S,M)$ yields a cluster-tilting object in $\mc_{(S,M)}$ and $\Ext^1_\mc(\gz,\gz)=0$ if $\gz$ is an internal arc in $(S,M).$
\end{cor}
\begin{pf}Let $\ggz$ be the triangulation of $(S,M)$ that we studied before, and $\ggz'$ be another triangulation. By Theorem \ref{fst}, there exists a sequence of flips which transform $\ggz$ to $\ggz',$
$$\ggz=\ggz_0\frac{f_{i_1}}{}\ggz_1\frac{f_{i_2}}{}\ggz_2\frac{f_{i_3}}{}\cdots\frac{f_{i_n}}{}\ggz_n=\ggz'$$
where $\ggz_j$ is a triangulation of $(S,M)$ and $f_{i_j}$ is a flip for $1\leq j\leq n.$
We know that $\ggz$ induces a cluster-tilting object $T_{\ggz}$ which is the direct sum of all internal arcs in $\ggz,$ and Theorem \ref{thm3} implies that $\ggz_1$ also induces a cluster-tilting object $T_{\ggz_1}$ given by all internal arcs in $\ggz_1$. The above lemma allows us to proceed to do this until $\ggz_n=\ggz'$ which completes the proof.
\end{pf}

Let $\ggz$ be a triangulation of $(S,M)$ and let $\gz$ and $\da$ be two non-boundary curves in $(S,M)$. From Theorem \ref{AKZ} we know that the morphism space $\Hom_\mc(\gz,\da)$ in $\mc=\mc_\ggz$ can be decomposed as
$k$-vector space as follows:

$$\Hom_{J(Q_\ggz,W_\ggz)}(M(\ggz,\gz),M(\ggz,\da)) \oplus  \{f\in \Hom_\mc(\gz,\da)|f \mbox{~~factors through~} T_\ggz\}.$$

Moreover, by Lemma 3.3 in \cite{Palu}, the morphisms factoring through $T_\ggz$ can be substituted by the morphisms of two other related objects in \mod $J(Q_\ggz,W_\ggz)$, so $\Hom_\mc(\gz,\da)$ is given by
$$\Hom_{J(Q_\ggz,W_\ggz)}(M(\ggz,\gz),M(\ggz,\da)) \oplus D\Hom_{J(Q_\ggz,W_\ggz)}(\t^{-1}M(\ggz,\da),\t M(\ggz,\gz))$$
In \cite{C-B} a basis of the Hom-space between string modules is explicitly described using factor strings and substrings. Using Theorem \ref{ABCP}, one can translate this into a description using curves and the triangulation $\ggz$ of $(S,M)$, thus one could give a basis of
  $\Hom_{\mc}(\gz,\da)$ combinatorially in terms of $\ggz$.

Another possibility is to choose a triangulation $\ggz'$ such that a given morphism in $\mc$ factoring through $T_\ggz$ can be described by a morphism in $\Hom_{J(Q_{\ggz'},W_{\ggz'})}(M(\ggz',\gz),M(\ggz',\da))$. This is illustrated in the following example.

\begin{example}\label{ex2}We consider the following two curves $\gz$, $\da$ from Example \ref{ex1},

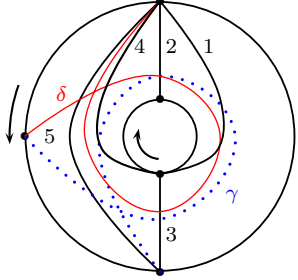
\begin{center}
\begin{pspicture}(-4,-1.8)(4,1.8)
\put(0,.5){\circle*{.1}}\put(0,-.5){\circle*{.1}}
\pscircle[linewidth=.7pt](0,0){.5}
\pscircle[linewidth=.7pt](0,0){1.8}
\put(0,1.8){\circle*{.1}}
\put(0,-1.8){\circle*{.1}}
\put(-1.8,0){\circle*{.1}}
\pscurve[linewidth=.7pt]{-}(0,1.8)(-1.2,0)(0,-1.8)
\uput[l](-1.2,0){$_{5}$}

\psarc[arcsepB=4pt]{<-}(0,0){.3}{-200}{-80}\psarc[arcsepB=4pt]{->}(0,0){2}{-200}{-175}
\psline[linewidth=.7pt]{-}(0,0.5)(0,1.8)\psline[linewidth=.7pt]{-}(0,-0.5)(0,-1.8)

\pscurve[linewidth=1.2pt,linestyle=dotted,linecolor=blue]{-}(-1.8,0)(0,-1.1)(1,0)(0,.8)(-.8,0)(0,-1.8)
\pscurve[linewidth=.5pt,linecolor=red]{-}(0,1.8)(-1,0)(0,-1)(.8,0)(0,.8)(-1.8,0)

\pscurve[linewidth=.7pt]{-}(0,1.8)(.8,0.5)(.75,-.2)(0,-.5)\pscurve{-}(0,1.8)(-.8,0.4)(-.75,-.2)(0,-.5)
\uput[r](.4,1.2){$_1$}{\red\uput[u](-1.3,.3){$_\da$}}
\uput[l](0,1.2){$_4$}
\uput[r](-.1,1.2){$_2$}
\uput[r](-.1,-1.3){$_3$}{\blue\uput[r](.7,-.8){$_\gz$}}
\end{pspicture}
\end{center}

where $\da$ is given by a full line and $\gz$ is given by a dotted line. It is easy to see that

\begin{center}
\begin{pspicture}(1.5,-.6)(6.3,0.5)
\put(1.5,0){$_{M(\ggz,\da)=}$}
\put(3.1,0.3){$_{1}$}
\put(2.8,0){$_{~3~2~,}$}
\put(3.2,-.3){$_{~4}$}\put(3.3,-.6){$_{~5}$}

\put(4.5,0){$_{M(\ggz,\gz)=}$}
\put(5.8,0.3){$_{5~1}$}
\put(5.9,0){$_{3~2~.}$}
\put(6.2,-.3){$_{~4}$}
\end{pspicture}
\end{center}

Hence $\Hom_{J(Q_\ggz,W_\ggz)}(M(\ggz,\gz),M(\ggz,\da))$ is one-dimensional with a basis element  $f$ mapping the factor string $5$ of $w(\ggz,\gz)$ to the substring $5$ of $w(\ggz,\da)$. From Corollary  \ref{boundary} one obtains that $\gz$ and $\da$ lie in a tube of rank 3 in  $\mc$ as follows:

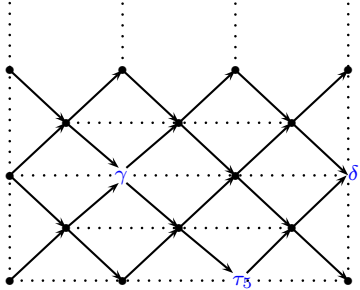
\begin{center}
\begin{pspicture}(-4,-.8)(4,3.1)
\put(-.75,-.7){\circle*{.1}}{\blue\put(.71,-.7){$_{\tau_5}$}}\put(-2.25,-.7){\circle*{.1}}\put(2.25,-.7){\circle*{.1}}
\put(.75,.7){\circle*{.1}}{\blue\put(-.85,.7){$_{\gz}$}}\put(-2.25,.7){\circle*{.1}}{\blue\put(2.25,.7){$_{\da}$}}
\put(1.5,0){\circle*{.1}}\put(-1.5,0){\circle*{.1}}\put(0,0){\circle*{.1}}
\put(1.5,1.4){\circle*{.1}}\put(-1.5,1.4){\circle*{.1}}\put(0,1.4){\circle*{.1}}
\put(-.75,2.1){\circle*{.1}}\put(.75,2.1){\circle*{.1}}\put(-2.25,2.1){\circle*{.1}}\put(2.25,2.1){\circle*{.1}}
\psline[linewidth=1pt,linestyle=dotted]{-}(2.25,-.7)(2.25,.6)\psline[linewidth=1pt,linestyle=dotted]{-}(2.25,3)(2.25,.9)
\psline[linewidth=1pt,linestyle=dotted]{-}(-2.25,-.7)(-2.25,3)
\psline{->}(-2.25,-.7)(-1.5,0)\psline{->}(-.75,-.7)(0,0)\psline{->}(.9,-.6)(1.5,0)
\psline{<-}(-.75,-.7)(-1.5,0)\psline{<-}(.72,-.6)(0,0)\psline{<-}(2.25,-.7)(1.5,0)

\psline{->}(-1.5,0)(-.8,0.6)\psline{->}(0,0)(.75,0.7)\psline{->}(1.5,0)(2.23,0.7)
\psline{->}(-2.25,.7)(-1.5,0)\psline{->}(-.7,.6)(0,0)\psline{->}(.75,.7)(1.5,0)

\psline{->}(-2.25,.7)(-1.5,1.4)\psline{->}(-.7,.8)(0,1.4)\psline{->}(.75,.7)(1.5,1.4)
\psline{->}(-1.5,1.4)(-.8,.8)\psline{->}(0,1.4)(.75,.7)\psline{->}(1.5,1.4)(2.25,.7)

\psline{->}(-1.5,1.4)(-.75,2.1)\psline{->}(0,1.4)(.75,2.1)\psline{->}(1.5,1.4)(2.25,2.1)
\psline{->}(-2.25,2.1)(-1.5,1.4)\psline{->}(-.75,2.1)(0,1.4)\psline{->}(.75,2.1)(1.5,1.4)
\psline[linewidth=1pt,linestyle=dotted]{-}(-.75,2.1)(-.75,3)\psline[linewidth=1pt,linestyle=dotted]{-}(.75,2.1)(.75,3)
\psline[linewidth=1pt,linestyle=dotted]{-}(-1.5,0)(1.5,0)\psline[linewidth=1pt,linestyle=dotted]{-}(-1.5,1.4)(1.5,1.4)
\psline[linewidth=1pt,linestyle=dotted]{-}(-2.25,.7)(-.9,0.7)\psline[linewidth=1pt,linestyle=dotted]{-}(2.1,.7)(-.6,0.7)
\psline[linewidth=1pt,linestyle=dotted]{-}(-2.25,-.7)(.6,-0.7)\psline[linewidth=1pt,linestyle=dotted]{-}(2.25,-.7)(.9,-0.7)
\end{pspicture}
\end{center}

Thus there is, besides the morphism $f$ induced from $\mod J(Q_\ggz,W_\ggz)$, another non-zero morphism $g\in \Hom_\mc(\gz,\da)$ factoring through $\tau_5$. Obviously, $g$ cannot be described by morphisms in $\mod J(Q_\ggz,W_\ggz)$ since $M(\ggz, \t_5)=0$, but we can realize $g$ in another module category: Applying a flip along $\tau_5$ we obtain a new triangulation $\ggz'$ whose associated quiver with potential is $Q_{\ggz'}=\tilde{A}_4$, $W_{\ggz'}=0$. Denote by $5^*=M(\ggz',\t_5)$ the simple $\tilde{A}_4-$module concentrated in 5, then

\begin{center}
\begin{pspicture}(-4.6,-.4)(3.4,0.5)

\put(-4.6,0){$_{M(\ggz',\da)=}$}
\put(-3.2,0.3){$_{~1}$}
\put(-3.35,0){$_{~3~2~,}$}
\put(-3.5,-.3){$_{~5^*}$} \put(-2.8,-.3){$_{4}$}

\put(-1.8,0){$_{M(\ggz',\gz)=}$}
\put(-.5,0.3){$_{~1}$}
\put(-.5,0){$_{3~2~5^*,}$}
\put(-.2,-.3){$_{~4}$}

\put(1.6,0){$_{M(\ggz',\tau_5)=}$}
\put(3,-.05){$_{5^*.}$}
\end{pspicture}
\end{center}

Hence $\Hom_{\tilde{A}_4}(M(\ggz',\gz),M(\ggz',{\da}))$ is also one dimensional given by a  basis element $g'$ mapping the factor string $5^*$ of $w(\ggz',\gz)$ to the substring $5^*$ of $w(\ggz',\da)$, which factors through $M(\ggz',\tau_5)$. Therefore $g\in \Hom_\mc(\gz,\da)$ can be described by $g'\in \Hom_{\tilde{A}_4}(M(\ggz',\gz),M(\ggz',\da))$. Moreover, since $\tilde{A}_4$ is hereditary, we can conclude that
\[
\begin{array}{rcl}
\Hom_\mc(\gz,\da) &\cong &\Hom_{\mc(\tilde{A}_4)}(M(\ggz',\gz),M(\ggz',\da))\\
 & =& \Hom_{\tilde{A}_4}(M(\ggz',\gz),M(\ggz',{\da}))\oplus \Ext^1_{\tilde{A}_4}(M(\ggz',\gz),M(\ggz',\gz))
\end{array}
\]
is two-dimensional. Hence
$$\Hom_\mc(\gz,\da) \cong \Hom_{J(Q_\ggz,W_\ggz)}(M(\ggz,\gz),M(\ggz,\da))\oplus \Hom_{\tilde{A}_4}(M(\ggz',\gz),M(\ggz',\da))$$
as $k$-vector space.
\end{example}
Using Lemma 3.3 in \cite{Palu}, we have \begin{center}
\begin{pspicture}(1.5,-.6)(6.3,0.5)
\put(-.2,0){$_{\t M(\ggz,\gz)=M(\ggz,\da)=}$}
\put(2.6,0.3){$_{1}$}
\put(2.3,0){$_{~3~2~,}$}
\put(2.7,-.3){$_{~4}$}\put(2.8,-.6){$_{~5}$}

\put(3.4,0){$_{\t^{-1}M(\ggz,\da)=M(\ggz,\gz)=}$}
\put(6.3,0.3){$_{5~1}$}
\put(6.4,0){$_{3~2~.}$}
\put(6.7,-.3){$_{~4}$}
\end{pspicture}
\end{center}
Hence, $\Hom_\mc(\gz,\da)$ is given by $$ \Hom_{J(Q_\ggz,W_\ggz)}(M(\ggz,\gz),M(\ggz,\da))\oplus D\Hom_{J(Q_\ggz,W_\ggz)}(M(\ggz,\gz),M(\ggz,\da))$$
as $k$-vector space.
\section{Extensions and Intersections}

We study in this section the relation between extensions in the category
$\mc=\mc_{(S,M)}$ and intersections of curves in $(S,M)$.
Given any curve in $(S,M)$ with self-intersections, we explicitly construct one
or two new curves, sometimes resolving the self-intersection, and sometimes
increasing the winding number.
These  curves serve as middle term of certain non-split short exact sequences
which allow to prove the following theorem:

\begin{thm}\label{intersection}$\Ext^1_\mc(\gz,\gz)=0$ if and only if $\gz$ is an internal arc in $(S,M).$
\end{thm}

\begin{pf}
Let  $\gz$ be a curve in $(S,M)$ and fix a triangulation $\ggz$ of $(S,M)$.
We denote by $\t_{i_1}, \t_{i_2},\ldots, \t_{i_d}$ the internal arcs of  $\ggz$ that intersect $\gz$ as indicated in  Figure 8:

\begin{center}
\begin{pspicture}(0,-2.1)(9.7,1.9)
\put(0,0){\circle*{.1}}\put(9.5,0){\circle*{.1}}
\put(1,1.5){\circle*{.1}}\put(2.5,1.5){\circle*{.1}}\put(4,1.5){\circle*{.1}}\put(5.5,1.5){\circle*{.1}}
\put(7,1.5){\circle*{.1}}\put(8.5,1.5){\circle*{.1}}
\put(1,-1.5){\circle*{.1}}\put(2.5,-1.5){\circle*{.1}}\put(3.5,-1.5){\circle*{.1}}\put(5,-1.5){\circle*{.1}}
\put(6.5,-1.5){\circle*{.1}}\put(8.5,-1.5){\circle*{.1}}
\psline[linewidth=1pt,linestyle=dotted]{-}(1.9,.3)(2.6,.3)\psline[linewidth=1pt,linestyle=dotted]{-}(4.8,.3)(5.5,.3)
\psline[linewidth=1pt,linestyle=dotted]{-}(7.5,.3)(8.2,.3)
\psline[linewidth=.5pt]{-}(2.8,-1.5)(2.5,-1.5)(1,-1.5)(0,0)(1,1.5)(1.5,1.5)
\psline[linewidth=.5pt]{-}(2,1.5)(2.5,1.5)(4,1.5)(4.5,1.5)
\psline[linewidth=.5pt]{-}(5,1.5)(5.5,1.5)(7,1.5)(7.5,1.5)
\psline[linewidth=.5pt]{-}(3,-1.5)(3.5,-1.5)(5,-1.5)(5.5,-1.5)
\psline[linewidth=.5pt]{-}(6,-1.5)(6.5,-1.5)(7.2,-1.5)
\psline[linewidth=.5pt]{-}(8,-1.5)(8.5,-1.5)(9.5,0)(8.5,1.5)(8,1.5)
\psline[linewidth=.5pt]{-}(1,-1.5)(1,1.5)(2.5,-1.5)\psline[linewidth=.5pt]{-}(2.5,1.5)(3.5,-1.5)(4,1.5)(5,-1.5)
\psline[linewidth=.5pt]{-}(5.5,1.5)(6.5,-1.5)(7,1.5)
\psline[linewidth=.5pt]{-}(8.5,1.5)(8.5,-1.5)
\psline[linewidth=1pt,linestyle=dotted]{-}(1.5,1.5)(2,1.5)\psline[linewidth=1pt,linestyle=dotted]{-}(4.5,1.5)(5,1.5)
\psline[linewidth=1pt,linestyle=dotted]{-}(7.5,1.5)(8,1.5)\psline[linewidth=1pt,linestyle=dotted]{-}(2.6,-1.5)(3,-1.5)
\psline[linewidth=1pt,linestyle=dotted]{-}(5.5,-1.5)(6,-1.5)\psline[linewidth=1pt,linestyle=dotted]{-}(7,-1.5)(8,-1.5)
\psline[linecolor=red,linewidth=.5pt](0,0)(9.5,0)

{\small\uput[l](0,0){$_{s(\gz)}$}\uput[r](9.5,0){$_{e(\gz)}$}
\uput[d](1.75,-1.4){$_{\t_{k_1}}$}\uput[d](4.25,-1.4){$_{\t_{k_{s+1}}}$}\uput[l](0.6,-.75){$_{\t_{k_0}}$}\uput[l](0.6,.75){$_{\t_{k_{-1}}}$}
\uput[u](3.25,1.4){$_{\t_{k_s}}$}\uput[u](6.25,1.4){$_{\t\scriptscriptstyle_{k_{s+n}}}$}\uput[r](9,.75){$_{\t_{k_d}}$}\uput[r](9,-.75){$_{\t_{k_{d+1}}}$}
\uput[l](1.2,.75){$_{\t_{i_1}}$}\uput[r](1.2,.75){$_{\t_{i_2}}$}\uput[l](2.75,.75){$_{\t_{i_s}}$}\uput[l](4,.75){$_{\t_{i_{s+1}}}$}
\uput[l](5.8,.75){$_{\t_{i_{s+n}}}$}\uput[r](6.8,.75){$_{\t_{i_{s+n+1}}}$}\uput[r](8.4,.75){$_{\t_{i_d}}$}\uput[r](4.2,.5){$_{\t_{i_{s+2}}}$}
\put(.6,-.5){$_{\vartriangle_0}$}\put(1.5,-.8){$_{\vartriangle_1}$}\put(3.25,.2){$_{\vartriangle_s}$}\put(4,-.8){$_{\vartriangle_{s+1}}$}
\put(6.1,.3){$_{\vartriangle_{s+n}}$}\put(8.7,-.5){$_{\vartriangle_d}$}}
{\red \uput[d](5.5,0){$_{\gz}$}}

\uput[d](5,-1.8){$_{Figure~8}$}
\end{pspicture}
\end{center}

Along its way, the curve $\gz$ is passing through (not necessarily distinct) triangles $ \vartriangle_0, \vartriangle_1, \ldots, \vartriangle_d$. For $1\leq l \leq d-1$, the triangle $\vartriangle_l$ is formed by the internal arcs $\t_{i_l}$ and $\t_{i_{l+1}}$ and a third arc which is denoted by $\t_{k_l}$.
In the first triangle $\vartriangle_0$, we denote the arcs clockwise by $\t_{i_1}$, $\t_{k_0}$ and $\t_{k_{-1}}$. Similarly,  $\vartriangle_d$ is formed by $\t_{i_d}$, $\t_{k_d}$ and $\t_{k_{d+1}}$, as shown in Figure~8.

By corollary \ref{cor1} we only need to prove one direction of the theorem, namely that $\Ext^1_\mc(\gz,\gz)\neq 0$ in case the curve $\gz$ has self-intersections. To do so, it is enough to prove $\Ext^1_{J(Q_{\ggz'},W_{\ggz'})}(M(\ggz',\gz),M(\ggz',\gz))\neq 0$ for some triangulation $\ggz'$. In some cases one can choose $\ggz'=\ggz$, but  sometimes it will be necessary to change $\ggz$ to $\ggz'$ by a sequence of flips in order to realize the extension over some algebra $J(Q_{\ggz'},W_{\ggz'})$.

We therefore assume that the curve $\gz$ has self-intersections, thus there are $0 < r < r' < 1$ such that $\gz(r)=\gz(r')$. We choose $r$ and $r'$ in such a way that the restriction $b=\gz|_{[ r,r']}$ is a simple closed curve. Hence the subword $w(\ggz,\gz|_{[ r,r']})$ of $w(\ggz,\gz)$ defines a band $w(\ggz,b)$ in $\mod J(Q_\ggz,W_\ggz)$ as follows:
\begin{center}
\begin{pspicture}(-3.8,-.5)(12,0.5)
\put(-.9,0){$_{w(\ggz,b):}$}
{\small\put(.2,0.1){$_{\t_{i_{s+1}}}$} \put(.9,-.5){$_{\t_{i_{s+n}}}$}
\put(.9,.5){$_{\t_{i_{s+2}}}$}
\uput[d](.6,-.1){$_{\theta}$}}
\psline{-}(.5,0.1)(.9,.4)
\psline{->}(.5,-.1)(.9,-.4)
\psline{-}(1.6,.5)(2.5,.5)
\psline[linestyle=dotted,linewidth=1pt](1.6,-.5)(2.1,-.5)
\psline{-}(2.2,-.5)(2.55,-.5)
\psline{-}(3.2,.3)(3.5,.1)
\psline{-}(3.1,-.4)(3.5,-.1)
\put(2.6,.5){$_{\t_{i_{s+3}}}$} \put(2.6,-.5){$_{\t_{i_{s+5}}}$}
\put(3.5,0){$_{\t_{i_{s+4}}}$}
\end{pspicture}
\end{center}
The difference between the string $w(\ggz,\gz|_{[ r,r']})$ and the band $w(\ggz,b)$ is the arrow between the endpoints. Up to duality, we may assume that this arrow $\theta$ (which is induced by the triangle $\vartriangle_{i_s}$), is oriented from $\tau_{i_{s+1}}$ towards $\tau_{i_{s+n}}$. We distinguish in the following several cases how the band  $w(\ggz,b)$ is embedded in the string $w(\ggz,\gz)$:
\medskip

{\bf Case I:} $w(\ggz,\gz)$ contains the band $w(\ggz,b).$

In this case, $w(\ggz,\gz)$ contains the string $w(\ggz,\gz|_{[ r,r']})$ and  the arrow $\theta$ and thus  $m \ge 1$ consecutive copies of $w(\ggz,b)$. We extend the subword of $w(\ggz,\gz)$ which is induced by $b$ to maximal length, that is, we choose $j \le s$ and $l \ge s+n$ such that all arrows between $\tau_{i_{j+1}}$ and $\tau_{i_l}$ belong to $w(\ggz,b)$, but the arrows $\tau_{i_{j}}\frac{\ \az\ }{}\tau_{i_{j+1}}$ induced by $\vartriangle_{i_j}$ and  $\tau_{i_{l}}\frac{\ \bz\ }{}\tau_{i_{l+1}}$ induced by $\vartriangle_{i_l}$ are not induced by $b$.

Without loss of generality, assume $\az$ is oriented as $\tau_{i_{j}}\overset{\az}{\ra}\tau_{i_{j+1}}$. Then $\bz$ has orientation $\tau_{i_{l}}\overset{\bz}{\ra}\tau_{i_{l+1}}$ since $\gz$ intersects itself. Therefore the arrows $\az':\tau_{i_{j+1}}\ra\t_{k_j}$ induced by $\vartriangle_{i_j}$ and the arrow $\bz':\t_{k_l}\ra\tau_{i_{l}}$ induced by $\vartriangle_{i_l}$ belong to $w(\ggz,b)$. We assume $\t_{i_{j+t}}=\t_{i_l},\t_{i_{j+t+1}}=\t_{k_l}$ with
$1\leq t\leq n.$ If $m=1$, then the situation is depicted as follows:

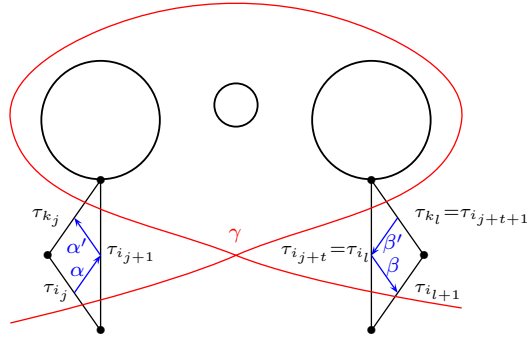
\begin{center}
\begin{pspicture}(-3,-2)(3,2.3)
    \pscircle[linewidth=.7pt](-1.8,0.8){.8}
    \pscircle[linewidth=.7pt](1.8,0.8){.8}
\put(-1.8,0){\circle*{.1}}\put(-1.8,-2){\circle*{.1}}\put(-2.5,-1){\circle*{.1}}
\put(1.8,0){\circle*{.1}}\put(1.8,-2){\circle*{.1}}\put(2.5,-1){\circle*{.1}}
\psline[linewidth=.5pt]{-}(-1.8,0)(-1.8,-2)(-2.5,-1)(-1.8,0)\psline[linewidth=.5pt]{-}(1.8,0)(1.8,-2)(2.5,-1)(1.8,0)
\psline[linewidth=.5pt,linecolor=blue]{->}(-1.8,-1)(-2.15,-.5)\psline[linewidth=.5pt,linecolor=blue]{->}(-2.15,-1.5)(-1.8,-1)
\psline[linewidth=.5pt,linecolor=blue]{<-}(1.8,-1)(2.15,-.5)\psline[linewidth=.5pt,linecolor=blue]{<-}(2.15,-1.5)(1.8,-1)
{\blue\uput[u](-2.1,-1.45){$_{\az}$}\uput[u](-2.1,-1.15){$_{\az'}$}
\uput[u](2.1,-1.45){$_{\bz}$}\uput[u](2.1,-1.15){$_{\bz'}$}}
\pscurve[linewidth=.5pt,linecolor=red]{-}(-3,-1.9)(0,-1)(3,0.8)(1.8,2)(-1.8,2)(-3,0.8)(0,-1)(3,-1.7)
{\red\uput[u](0,-1){$_{\gz}$}}
{\small\uput[r](-1.9,-1){$_{\t_{i_{j+1}}}$}\uput[l](-2,-1.5){$_{\t_{i_j}}$}\uput[u](-2.5,-.8){$_{\t_{k_j}}$}
\uput[l](2,-1){$_{\t_{i_{j+t}}=\t_{i_l}}$}\uput[r](2.2,-1.5){$_{\t_{i_{l+1}}}$}\uput[r](2.2,-.5){$_{\t_{k_{l}}=\t_{i_{j+t+1}}}$}}
\pscircle[linewidth=.7pt](0,1){.3}
\end{pspicture}
\end{center}
Here, the boundaries are indicated by the circles, and $\t_{i_j}=\t_{k_0}$ when $j=0$ (also note that $\t_{i_j}$ might be a boundary arc). We rewrite Figure 8  accordingly:

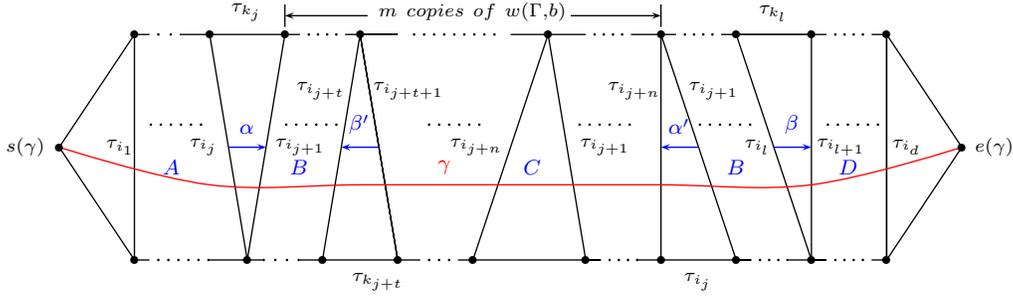
\begin{center}
\begin{pspicture}(-.2,-2)(14,2)

\psline[linewidth=.5pt]{-}(1.2,-1.5)(1,-1.5)(0,0)(1,1.5)(1.2,1.5)
\psline[linewidth=.5pt]{-}(1.6,1.5)(2,1.5)(3,1.5)(2.5,-1.5)(2,1.5)\psline[linewidth=.5pt]{-}(3.2,1.5)(3,1.5)
\psline[linewidth=.5pt]{-}(2,-1.5)(2.7,-1.5)\psline[linewidth=.5pt]{-}(3.2,-1.5)(4.8,-1.5)
\psline[linewidth=.5pt]{-}(3.5,-1.5)(4,1.5)(4.5,-1.5)(4,1.5)(4.5,1.5)(3.8,1.5)
\psline[linewidth=.5pt]{-}(5.3,-1.5)(7.2,-1.5)(7,-1.5)(6.5,1.5)(5.5,-1.5)\psline[linewidth=.5pt]{-}(6,1.5)(7,1.5)
\psline[linewidth=.5pt]{-}(7.7,-1.5)(9.2,-1.5)(9,-1.5)(8,1.5)(8,-1.5)\psline[linewidth=.5pt]{-}(7.6,1.5)(8.2,1.5)
\psline[linewidth=.5pt]{-}(9.7,-1.5)(10.2,-1.5)(10,-1.5)(10,1.5)(9,1.5)(10,-1.5)
\psline[linewidth=.5pt]{-}(10.8,-1.5)(11,-1.5)(12,0)(11,1.5)(11,-1.5)(11,1.5)(10.8,1.5)
\psline[linewidth=.5pt]{-}(8.8,1.5)(10.2,1.5)\psline[linewidth=.5pt]{-}(1,1.5)(1,-1.5)
\psline[linewidth=.5pt]{-}(3,1.6)(3,1.9)\psline[linewidth=.5pt]{-}(8,1.6)(8,1.9)
\psline[linewidth=.5pt]{<-}(3,1.75)(4.1,1.75)\psline[linewidth=.5pt]{->}(6.8,1.75)(8,1.75)
{\blue\uput[d](1.5,0){$_{A}$}\uput[d](3.2,0){$_{B}$}\uput[d](6.3,0){$_{C}$}
\uput[d](9,0){$_{B}$}\uput[d](10.5,0){$_{D}$}}
\pscurve[linewidth=.5pt,linecolor=red]{-}(0,0)(2,-.5)(4,-.5)(8,-.5)(10,-.5)(12,0)
{\small\uput[u](5.5,1.5){$_{m~copies~ of~ w(\ggz,b)}$}
\uput[l](0,0){$_{s(\gz)}$}\uput[r](12,0){$_{e(\gz)}$}
\uput[u](2.5,1.5){$_{\t_{k_{j}}}$}
\uput[d](4.25,-1.5){$_{\t_{k_{j+t}}}$}
\uput[l](1.2,0){$_{\t_{i_1}}$}
\uput[l](2.3,0){$_{\t_{i_j}}$}\uput[r](2.7,0){$_{\t_{i_{j+1}}}$}
\uput[d](8.5,-1.5){$_{\t_{i_j}}$}
\uput[l](4,.75){$_{\t_{i_{j+t}}}$}\uput[r](10.9,0){$_{\t_{i_d}}$}\uput[r](4,.75){$_{\t_{i_{j+t+1}}}$}
\uput[l](8.2,.75){$_{\t_{i_{j+n}}}$}\uput[r](8.2,.75){$_{\t_{i_{j+1}}}$}
\uput[l](9.65,0){$_{\t_{i_{l}}}$}
\uput[l](6.1,0){$_{\t_{i_{j+n}}}$}\uput[r](6.75,0){$_{\t_{i_{j+1}}}$}\uput[u](9.5,1.5){$_{\t_{k_l}}$}\uput[r](9.9,0){$_{\t_{i_{l+1}}}$}
}
\put(0,0){\circle*{.1}}\put(12,0){\circle*{.1}}
\put(1,1.5){\circle*{.1}}\put(2,1.5){\circle*{.1}}\put(3,1.5){\circle*{.1}}\put(4,1.5){\circle*{.1}}
\put(6.5,1.5){\circle*{.1}}\put(8,1.5){\circle*{.1}}\put(9,1.5){\circle*{.1}}\put(10,1.5){\circle*{.1}}
\put(11,1.5){\circle*{.1}}
\psline[linewidth=1pt,linestyle=dotted](1.3,1.5)(1.5,1.5)\psline[linewidth=1pt,linestyle=dotted](3.3,1.5)(3.7,1.5)
\psline[linewidth=1pt,linestyle=dotted](4.7,1.5)(5.8,1.5)\psline[linewidth=1pt,linestyle=dotted](7.1,1.5)(7.5,1.5)
\psline[linewidth=1pt,linestyle=dotted](8.3,1.5)(8.7,1.5)\psline[linewidth=1pt,linestyle=dotted](10.3,1.5)(10.7,1.5)
\put(1,-1.5){\circle*{.1}}\put(2.5,-1.5){\circle*{.1}}\put(3.5,-1.5){\circle*{.1}}\put(4.5,-1.5){\circle*{.1}}
\put(5.5,-1.5){\circle*{.1}}\put(7,-1.5){\circle*{.1}}\put(8,-1.5){\circle*{.1}}\put(9,-1.5){\circle*{.1}}
\put(10,-1.5){\circle*{.1}}\put(11,-1.5){\circle*{.1}}
\psline[linewidth=1pt,linestyle=dotted]{-}(1.3,-1.5)(2,-1.5)\psline[linewidth=1pt,linestyle=dotted]{-}(2.7,-1.5)(3.2,-1.5)
\psline[linewidth=1pt,linestyle=dotted]{-}(4.9,-1.5)(5.2,-1.5)\psline[linewidth=1pt,linestyle=dotted]{-}(7.3,-1.5)(7.7,-1.5)
\psline[linewidth=1pt,linestyle=dotted]{-}(9.7,-1.5)(9.3,-1.5)\psline[linewidth=1pt,linestyle=dotted]{-}(10.3,-1.5)(10.7,-1.5)
\psline[linewidth=1pt,linestyle=dotted]{-}(1.2,.3)(1.9,.3)
\psline[linewidth=1pt,linestyle=dotted]{-}(3,.3)(3.7,.3)
\psline[linewidth=1pt,linestyle=dotted]{-}(4.9,.3)(5.6,.3)
\psline[linewidth=1pt,linestyle=dotted]{-}(6.9,.3)(7.6,.3)
\psline[linewidth=1pt,linestyle=dotted]{-}(8.5,.3)(9.2,.3)
\psline[linewidth=1pt,linestyle=dotted]{-}(10.2,.3)(10.9,.3)
\psline[linewidth=.5pt,linecolor=blue]{->}(2.25,0)(2.75,0)\psline[linewidth=.5pt,linecolor=blue]{->}(4.25,0)(3.75,0)
\psline[linewidth=.5pt,linecolor=blue]{<-}(8,0)(8.5,0)\psline[linewidth=.5pt,linecolor=blue]{->}(9.5,0)(10,0)
{\blue\uput[u](2.5,0){$_{\az}$}\uput[u](4,0){$_{\bz'}$}\uput[u](8.25,0){$_{\az'}$}\uput[u](9.75,0){$_{\bz}$}}
{\red\uput[d](5,0){$_{\gz}$}}
\end{pspicture}
\end{center}
Thus $w(\ggz,\gz)$ is given by
{\small $$w(\ggz,\gz)=A-\t_{i_j}\overset{\az}{\ra}\t_{i_{j+1}}-B-\t_{i_l}\overset{\bz'}{\leftarrow} \t_{k_l}-C-\t_{i_{j+n}}\overset{\az'}{\leftarrow} \t_{i_{j+1}}-B-\t_{i_l}\overset{\bz}{\ra}\t_{i_{l+1}}-D.$$ }
Consider following two strings in $\mod J(Q_\ggz,Q_\ggz):$
{\small $$w(\ggz,\gz')=A-\t_{i_j}\overset{\az}{\ra}\t_{i_{j+1}}-B-\t_{i_l}\overset{\bz'}{\leftarrow} \t_{k_l}-C-\t_{i_{j+n}}\overset{\az'}{\leftarrow} \t_{i_{j+1}}-B-\t_{i_l}\overset{\bz'}{\leftarrow}\t_{k_l}$$
$$-C-\t_{i_{j+n}}\overset{\az'}{\leftarrow} \t_{i_{j+1}}-B-\t_{i_l}\overset{\bz}{\ra}\t_{i_{l+1}}-D.$$
$$w({\ggz,\gz''})=A-\t_{i_j}\overset{\az}{\ra}\t_{i_{j+1}}-B-\t_{i_l}\overset{\bz}{\ra}\t_{i_{l+1}}-D.$$}
where $\gz'$ and $\gz''$ are the corresponding curves in $(S,M)$. If $m=1$, then $\gz', \gz''$ can be visualized as follows where $\gz'$ is given by the full curve and  $\gz''$ is described by the dotted curve.
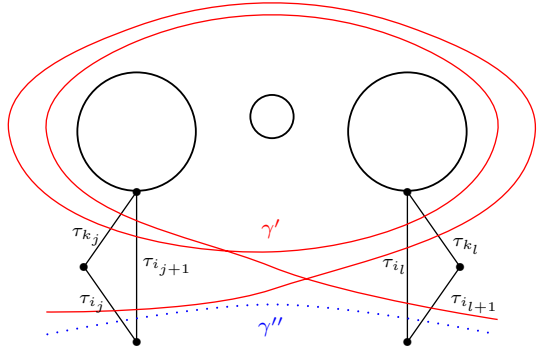
\begin{center}
\begin{pspicture}(-3,-2)(3,2.5)
\pscircle[linewidth=.7pt](-1.8,0.8){.8}    \pscircle[linewidth=.7pt](1.8,0.8){.8}
\pscircle[linewidth=.7pt](0,1){.3}

\put(-1.8,0){\circle*{.1}}\put(-1.8,-2){\circle*{.1}}\put(-2.5,-1){\circle*{.1}}
\put(1.8,0){\circle*{.1}}\put(1.8,-2){\circle*{.1}}\put(2.5,-1){\circle*{.1}}
\psline[linewidth=.5pt]{-}(-1.8,0)(-1.8,-2)(-2.5,-1)(-1.8,0)\psline[linewidth=.5pt]{-}(1.8,0)(1.8,-2)(2.5,-1)(1.8,0)
\pscurve[linewidth=.5pt,linecolor=red]{-}(-3,-1.6)(0,-1.3)(3.5,0.8)(1.8,2.2)(-1.8,2.2)(-3.5,.8)(0,-.8)(3,0.8)(1.8,2)(-1.8,2)(-3,0.8)(0,-1)(3,-1.7)
{\red\uput[u](0,-.8){$_{\gz'}$}}

{\small\uput[r](-1.9,-1){$_{\t_{i_{j+1}}}$}\uput[l](-2,-1.5){$_{\t_{i_j}}$}\uput[u](-2.45,-.9){$_{\t_{k_j}}$}
\uput[l](2,-1){$_{\t_{i_l}}$}\uput[r](2.2,-1.5){$_{\t_{i_{l+1}}}$}\uput[u](2.6,-1){$_{\t_{k_{l}}}$}}

\pscurve[linewidth=.8pt,linecolor=blue,linestyle=dotted]{-}(-3,-1.9)(0,-1.5)(3,-1.9)
{\blue\uput[d](0,-1.5){$_{\gz''}$}}
\end{pspicture}
\end{center}
We easily construct a non-split exact sequence in $\mod J(Q_\ggz,W_\ggz)$
$$0\lra M(\ggz,\gz)\overset{{\tiny _{\left[
                              \begin{array}{c}
                                f \\
                                g \\
                              \end{array}
                            \right]
}}}{\lra} M({\ggz,\gz'})\oplus M({\ggz,\gz''})\lra M(\ggz,\gz)\lra 0$$
where $f$ identifies the factor string
$$A-\t_{i_j}\overset{\az}{\ra}\t_{i_{j+1}}-B-\t_{i_l}\overset{\bz'}{\leftarrow} \t_{k_l}-C-\t_{i_{j+n}}\overset{\az'}{\leftarrow} \t_{i_{j+1}}-B-\t_{i_l}$$
 of $w(\ggz,\gz)$ with a  substring of $w(\ggz,\gz')$ and  $g$ sends the factor string \\$\t_{i_{j+1}}-B-\t_{i_l}\overset{\bz}{\ra}\t_{i_{l+1}}-D$ of $w(\ggz,\gz)$ to the same  substring of $w(\ggz,\gz'').$
\medskip

{\bf Case II:}  $w(\ggz,\gz)$ does not contain the band $w(\ggz,b).$

In this case, $w(\ggz,\gz)$ does not contain $\theta$ but it contains $w(\ggz,\gz|_{[ r,r']})$ as a subword. We distinguish a number of subcases:
\begin{itemize}
\item[\bf (II.1)]$0=s<s+n=d$.

In this subcase, $\t_{i_s}=\t_{k_0}$ and the subword $w(\ggz,\gz|_{[ r,r']})$ equals $w(\ggz,\gz)$, see the following picture
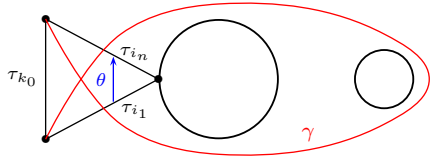
\begin{center}
\begin{pspicture}(0,-1.1)(9,1.1)
\pscircle[linewidth=.7pt](4.8,0){.8}
\pscircle[linewidth=.7pt](7,0){.4}
\put(4,0){\circle*{.1}}\put(2.5,.8){\circle*{.1}}
\put(2.5,-.8){\circle*{.1}}
\psline[linewidth=.5pt]{-}(2.5,.8)(4,0)(2.5,-.8)(2.5,.8)
\pscurve[linewidth=.5pt,linecolor=red]{-}(2.5,-.8)(3.5,.6)(4.8,1)(7,.6)(7.6,0)(7,-0.6)(4.8,-1)(3.5,-.6)(2.5,.8)
{\red\uput[u](6,-1){$_{\gz}$}}
{\small\uput[u](3.7,-.7){$_{\t_{i_{1}}}$}\uput[d](3.7,.6){$_{\t_{i_n}}$}
\uput[l](2.6,0){$_{\t_{k_0}}$}}
\psline[linewidth=.5pt,linecolor=blue]{->}(3.4,-.3)(3.4,0.3)
{\small\blue\uput[l](3.5,0){$_{\theta}$}}
\end{pspicture}
\end{center}
Since $w(\ggz,b)$ is a band we know that there exists a string of the following form  in $\mod J(Q_\ggz,W_\ggz)$
$$w(\ggz,\gz')=\t_{i_1}\leftrightsquigarrow \t_{i_n}\overset{\theta}{\leftarrow}  \t_{i_1}\leftrightsquigarrow \t_{i_n}$$
where $\gz'$ can be visualized as follows
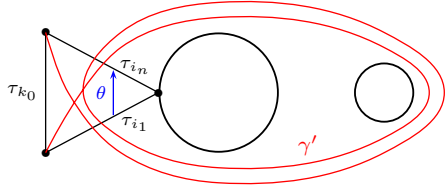
\begin{center}
\begin{pspicture}(0,-1.1)(9,1.1)
\pscircle[linewidth=.7pt](4.8,0){.8}
\pscircle[linewidth=.7pt](7,0){.4}
\put(4,0){\circle*{.1}}\put(2.5,.8){\circle*{.1}}
\put(2.5,-.8){\circle*{.1}}
\psline[linewidth=.5pt]{-}(2.5,.8)(4,0)(2.5,-.8)(2.5,.8)
\pscurve[linewidth=.5pt,linecolor=red,showpoints=false]{-}(2.5,-.8)(3.5,.6)(4.8,1)(7,.6)(7.6,0)(7,-0.6)(4.8,-1)(3.5,-.6)
(3,0)(3.5,.8)(4.8,1.2)(7,.8)(7.8,0)(7,-.8)(4.8,-1.2)(3.5,-.8)(2.8,0)(2.5,.8)
{\red\uput[u](6,-1){$_{\gz'}$}}
{\small\uput[u](3.7,-.7){$_{\t_{i_{1}}}$}\uput[d](3.7,.6){$_{\t_{i_n}}$}
\uput[l](2.6,0){$_{\t_{k_0}}$}}
\psline[linewidth=.5pt,linecolor=blue]{->}(3.4,-.3)(3.4,0.3)
{\small\blue\uput[l](3.5,0){$_{\theta}$}}
\end{pspicture}
\end{center}
It is easy to see that there is a non-split exact sequence in $\mod J(Q_\ggz,W_\ggz)$:
$$0\lra M(\ggz,\gz)\overset{f}{\lra} M(\ggz,\gz')\lra M(\ggz,\gz)\lra 0$$
where $f$ is induced by the embedding of the string $w(\ggz,\gz)$ as a substring of $w(\ggz,\gz')$.
\medskip

\item[\bf (II.2)]$0=s<s+n<d.$

In this subcase, $\t_{i_s}=\t_{i_0}=\t_{k_0}$ and the subword $w(\ggz,\gz|_{[ r,r']})$ is a proper subword of $w(\ggz,\gz)$ which ends before reaching the endpoint $e(\gz)$ of $\gz$, see the following figure

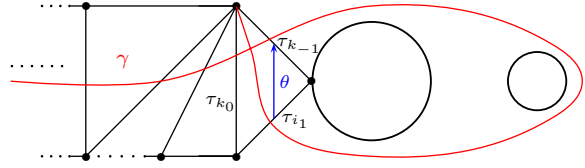
\begin{center}
\begin{pspicture}(0,-1.1)(9,1.1)
\pscircle[linewidth=.7pt](4.8,0){.8}    \pscircle[linewidth=.7pt](7,0){.4}
\put(4,0){\circle*{.1}}\put(3,1){\circle*{.1}}
\put(1,-1){\circle*{.1}}\put(1,1){\circle*{.1}}
\put(2,-1){\circle*{.1}}
\put(3,-1){\circle*{.1}}
\psline[linewidth=.5pt]{-}(2,-1)(3,1)(1,-1)(1,1)\psline[linewidth=.5pt]{-}(1,1)(.8,1)\psline[linewidth=.5pt]{-}(1,-1)(.8,-1)
\psline[linewidth=1pt,linestyle=dotted]{-}(.4,-1)(.8,-1)\psline[linewidth=1pt,linestyle=dotted]{-}(.4,1)(.8,1)
\psline[linewidth=.5pt]{-}(2.5,1)(3,1)(4,0)(3,-1)(2.5,-1)(3,-1)(3,1)
\psline[linewidth=1pt,linestyle=dotted]{-}(0,.2)(.7,.2)

\pscurve[linewidth=.5pt,linecolor=red]{-}(0,0)(2,0)(4.8,1)(7,.6)(7.6,0)(7,-0.6)(4.8,-1)(3.5,-.6)(3.2,.4)(3,1)
{\red\uput[u](1.5,0){$_{\gz}$}}
\psline[linewidth=.5pt]{-}(1,1)(3,1)
\psline[linewidth=.5pt]{-}(1.8,-1)(3,-1)
\psline[linewidth=1pt,linestyle=dotted]{-}(1,-1)(1.8,-1)
{\small\uput[u](3.8,-.8){$_{\t_{i_{1}}}$}\uput[d](3.85,.75){$_{\t_{k_{-1}}}$}
\uput[l](3.2,-.3){$_{\t_{k_0}}$}}
\psline[linewidth=.5pt,linecolor=blue]{->}(3.5,-.5)(3.5,0.5)
{\small\blue\uput[r](3.4,0){$_{\theta}$}}
\end{pspicture}
\end{center}
Since $\gz$ intersects itself, $\t_{k_{-1}}=\t_{i_n}$ and $\t_{i_{n+1}}=\t_{k_0}$. We rewrite Figure 8  accordingly:

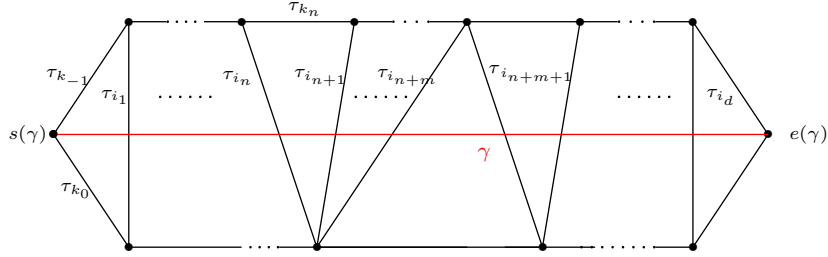
\begin{center}
\begin{pspicture}(0,-1.5)(9.7,1.7)
\put(0,0){\circle*{.1}}\put(9.5,0){\circle*{.1}}
\put(1,1.5){\circle*{.1}}\put(2.5,1.5){\circle*{.1}}\put(4,1.5){\circle*{.1}}\put(5.5,1.5){\circle*{.1}}
\put(7,1.5){\circle*{.1}}\put(8.5,1.5){\circle*{.1}}
\psline[linewidth=.5pt]{-}(2.5,-1.5)(1,-1.5)(0,0)(1,1.5)(1.5,1.5)
\psline[linewidth=.5pt]{-}(2,1.5)(2.5,1.5)(4,1.5)(4.5,1.5)
\psline[linewidth=.5pt]{-}(5,1.5)(5.5,1.5)(7,1.5)(7.5,1.5)
\psline[linewidth=.5pt]{-}(3,-1.5)(3.5,-1.5)(5.5,-1.5)
\psline[linewidth=.5pt]{-}(6,-1.5)(6.5,-1.5)(7.2,-1.5)
\psline[linewidth=.5pt]{-}(8,-1.5)(8.5,-1.5)(9.5,0)(8.5,1.5)(8,1.5)
\psline[linewidth=.5pt]{-}(1,-1.5)(1,1.5)\psline[linewidth=.5pt]{-}(2.5,1.5)(3.5,-1.5)(4,1.5)
\psline[linewidth=.5pt]{-}(5.5,1.5)(6.5,-1.5)(7,1.5)
\psline[linewidth=.5pt]{-}(8.5,1.5)(8.5,-1.5)
\psline[linewidth=.5pt]{-}(3.5,-1.5)(5.5,1.5)\psline[linewidth=.5pt]{-}(3.5,-1.5)(6.5,-1.5)
\psline[linewidth=1pt,linestyle=dotted]{-}(4,.5)(4.7,.5)
\psline[linewidth=1pt,linestyle=dotted]{-}(1.4,.5)(2.1,.5)
\psline[linewidth=1pt,linestyle=dotted]{-}(7.5,.5)(8.2,.5)
\psline[linewidth=1pt,linestyle=dotted]{-}(1.5,1.5)(1.9,1.5)
\psline[linewidth=1pt,linestyle=dotted]{-}(4.5,1.5)(5,1.5)
\psline[linewidth=1pt,linestyle=dotted]{-}(7.5,1.5)(8,1.5)\psline[linewidth=1pt,linestyle=dotted]{-}(2.6,-1.5)(3,-1.5)
\psline[linewidth=1pt,linestyle=dotted]{-}(7,-1.5)(8,-1.5)
\psline[linecolor=red,linewidth=.5pt](0,0)(9.5,0)

\put(1,-1.5){\circle*{.1}}{\circle*{.1}}\put(3.4,-1.5){\circle*{.1}}
\put(6.4,-1.5){\circle*{.1}}\put(8.4,-1.5){\circle*{.1}}

{\small\uput[l](0,0){$_{s(\gz)}$}\uput[r](9.5,0){$_{e(\gz)}$}
\uput[l](0.6,-.75){$_{\t_{k_0}}$}\uput[l](0.6,.75){$_{\t_{k_{-1}}}$}
\uput[u](3.25,1.4){$_{\t_{k_n}}$}
\uput[l](1.1,.5){$_{\t_{i_1}}$}\uput[l](2.75,.75){$_{\t_{i_n}}$}\uput[l](4,.75){$_{\t_{i_{n+1}}}$}
\uput[l](5.2,.75){$_{\t_{i_{n+m}}}$}\uput[r](5.5,.8){$_{\t\scriptscriptstyle_{i_{n+m+1}}}$}\uput[r](8.4,.5){$_{\t_{i_d}}$}}
{\red \uput[d](5.5,0){$_{\gz}$}}

\end{pspicture}
\end{center}

where $\t_{i_{n+2}},\t_{i_{n+3}}\ldots \t_{i_{n+m}}$ with $m \geq 1$ are all internal arcs lying clockwise before $\t_{i_{n+1}}$. After applying flips along $\t_{i_{n+1}},\t_{i_{n+2}}\ldots\t_{i_{n+m}}$, we get a new triangulation $\ggz'$ and a new band $w(\ggz',b)$ in $\mod J(Q_{\ggz'},W_{\ggz'})$ related to same closed curve $b$:

\begin{center}
\begin{pspicture}(-3.4,-.6)(12,0.6)
\put(-1,0){$_{w(\ggz',b):}$}
\put(.2,0.1){$_{\t^*_{i_{n+m}}}$}
\put(1,-.5){$_{\t_{i_n}}$}
\put(3.5,.1){$_{\t_{i_1}.}$}
\psline{-}(3.4,0.1)(3,.26)
\put(.9,.5){$_{\t^*_{i_{n+m-1}}}$}
\psline{-}(.5,0.1)(.9,.4)
\psline{-}(.6,-.2)(.9,-.4)
\psline[linestyle=dotted,linewidth=1pt](1.6,-.5)(2.6,-.5)
\psline[linestyle=dotted,linewidth=1pt](1.6,.5)(2.3,.5)
\psline{-}(2.8,-.5)(3.4,0)
\put(2.5,.5){$_{\t^*_{i_{n+1}}}$}
\end{pspicture}
\end{center}

If $d=n+m$, then the subword $w(\ggz',b)$ is equal to $w(\ggz',\gz)=\t^*_{i_{n+m}}-\t^*_{i_{n+m-1}}\leftrightsquigarrow \t_{i_n}$ in $\mod J(Q_{\ggz'},W_{\ggz'})$, and this is the same case as in II.1 when we consider $\mod J(Q_{\ggz'},W_{\ggz'})$. If $d>n+m$, then
$$w(\ggz',\gz)=\t^*_{i_{n+m}}-\t^*_{i_{n+m-1}}\leftrightsquigarrow \t_{i_n}- \t^*_{i_{n+m}}-\cdots,$$
hence $w(\ggz',\gz)$ contains the band $w(\ggz',b)$ as subword, this is the same case as in I.

\item[\bf (II.3)] $0<s<s+n<d.$

In this subcase, the subword $w(\ggz,\gz|_{[ r,r']})$ is a proper subword of $w(\ggz,\gz)$ which starts after $\t_{i_1}$ and ends before reaching the endpoint $e(\gz)$ of $\gz$. Since $\gz$ intersects itself, there must exist $1\leq t\leq s$ such that $\t_{i_{s+n+l}}=\t_{i_{s+1-l}}$ for $1 \leq l\leq s-t+1$ and $\t_{i_{s+n+t+1}}=\t_{k_{t-1}},$ see the following figure

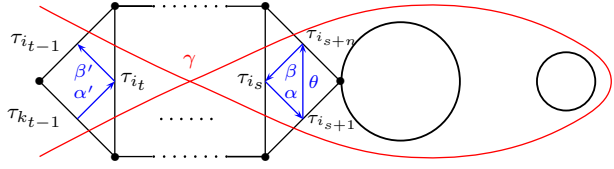
\begin{center}
\begin{pspicture}(0,-1.1)(9,1.1)
\pscircle[linewidth=.7pt](4.8,0){.8}    \pscircle[linewidth=.7pt](7,0){.4}
\put(4,0){\circle*{.1}}\put(3,1){\circle*{.1}}\put(1,1){\circle*{.1}}
\put(0,0){\circle*{.1}}\put(1,-1){\circle*{.1}}\put(3,-1){\circle*{.1}}

\psline[linewidth=1pt,linestyle=dotted]{-}(1.6,-.5)(2.3,-.5)

\psline[linewidth=.5pt]{-}(1.5,-1)(1,-1)(0,0)(1,1)(1.5,1)(1,1)(1,-1)\psline[linewidth=.5pt]{-}(2.5,1)(3,1)(4,0)(3,-1)(2.5,-1)(3,-1)(3,1)
\pscurve[linewidth=.5pt,linecolor=red]{-}(0,-1)(2,0)(4.8,1)(7,.6)(7.6,0)(7,-0.6)(4.8,-1)(2,0)(0,1)
{\red\uput[u](2,0){$_{\gz}$}}
{\small\uput[u](3.9,-.85){$_{\t_{i_{s+1}}}$}\uput[d](3.9,.85){$_{\t_{i_{s+n}}}$}}
\uput[r](.9,0){$_{\t_{i_t}}$}
\psline[linewidth=1pt,linestyle=dotted]{-}(1.5,1)(2.5,1)\psline[linewidth=1pt,linestyle=dotted]{-}(1.5,-1)(2.5,-1)
\uput[l](.5,.5){$_{\t_{i_{t-1}}}$}\uput[l](.5,-.5){$_{\t_{k_{t-1}}}$}
\uput[l](3.2,0){$_{\t_{i_s}}$}
\psline[linewidth=.5pt,linecolor=blue]{->}(3.5,.5)(3,0)
\psline[linewidth=.5pt,linecolor=blue]{<-}(3.5,-.5)(3,0)
\psline[linewidth=.5pt,linecolor=blue]{->}(3.5,-.5)(3.5,0.5)
{\small\blue\uput[r](3.4,0){$_{\theta}$}\uput[u](3.35,-.4){$_{\az}$}\uput[d](3.35,.4){$_{\bz}$}}
\psline[linewidth=.5pt,linecolor=blue]{->}(.5,-.5)(1,0)
\psline[linewidth=.5pt,linecolor=blue]{->}(1,0)(.5,0.5)
{\small\blue\uput[u](.6,-.4){$_{\az'}$}\uput[d](.6,.4){$_{\bz'}$}}
\end{pspicture}
\end{center}

where $\az$ and $\bz$ are arrows induced by the triangle $\vartriangle_{i_s}$, and $\az'$ and $\bz'$ are induced by the triangle $\vartriangle_{i_{t-1}}$. We rewrite Figure 8  accordingly:

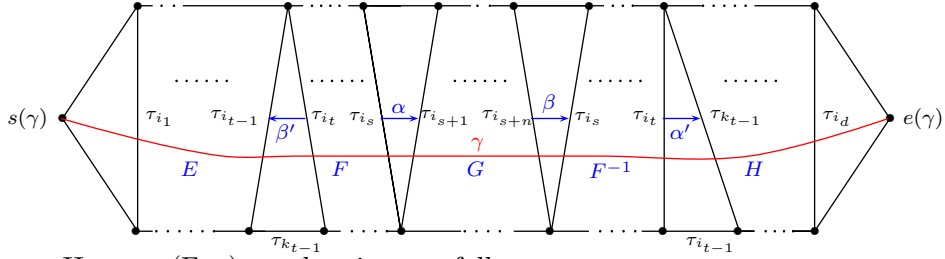
\begin{center}
\begin{pspicture}(0,-1.6)(14,1.5)
\put(0,0){\circle*{.1}}\put(11,0){\circle*{.1}}
\put(1,1.5){\circle*{.1}}\put(3,1.5){\circle*{.1}}\put(4,1.5){\circle*{.1}}\put(5,1.5){\circle*{.1}}\put(10,-1.5){\circle*{.1}}
\put(6,1.5){\circle*{.1}}\put(7,1.5){\circle*{.1}}\put(8,1.5){\circle*{.1}}\put(10,1.5){\circle*{.1}}
\psline[linewidth=1pt,linestyle=dotted](1.3,1.5)(1.5,1.5)\psline[linewidth=1pt,linestyle=dotted](3.3,1.5)(3.7,1.5)
\psline[linewidth=1pt,linestyle=dotted](5.5,1.5)(5.7,1.5)\psline[linewidth=1pt,linestyle=dotted](7.3,1.5)(7.6,1.5)
\psline[linewidth=1pt,linestyle=dotted](9.1,1.5)(8.5,1.5)
\psline[linewidth=.5pt]{-}(1.2,-1.5)(1,-1.5)(0,0)(1,1.5)(1.2,1.5)
\psline[linewidth=.5pt]{-}(1.6,1.5)(3,1.5)(2.5,-1.5)\psline[linewidth=.5pt]{-}(3.2,1.5)(3,1.5)
\psline[linewidth=.5pt]{-}(2,-1.5)(2.7,-1.5)
\psline[linewidth=.5pt]{-}(4,1.5)(4.5,-1.5)(4,1.5)(4.5,1.5)(3.8,1.5)
\psline[linewidth=.5pt]{-}(3.5,-1.5)(3,1.5)\psline[linewidth=.5pt]{-}(4.5,-1.5)(5,1.5)
\psline[linewidth=.5pt]{-}(6,1.5)(6.5,-1.5)(7,1.5)
\psline[linewidth=.5pt]{-}(6,1.5)(7,1.5)
\psline[linewidth=.5pt]{-}(7.7,-1.5)(9.2,-1.5)(9,-1.5)(8,1.5)(8,-1.5)\psline[linewidth=.5pt]{-}(7.6,1.5)(8.5,1.5)
\psline[linewidth=.5pt]{-}(9.2,1.5)(10,1.5)(11,0)(10,-1.5)(9.7,-1.5)(10,-1.5)(10,1.5)
\psline[linewidth=.5pt]{-}(1,1.5)(1,-1.5)\psline[linewidth=.5pt]{-}(6,-1.5)(7,-1.5)\psline[linewidth=.5pt]{-}(5.7,1.5)(7.2,1.5)
\psline[linewidth=.5pt]{-}(2.5,-1.5)(3.7,-1.5)
\psline[linewidth=.5pt]{-}(4.3,-1.5)(5.3,-1.5)
\uput[l](0,0){$_{s(\gz)}$}\uput[r](11,0){$_{e(\gz)}$}
\psline[linewidth=1pt,linestyle=dotted]{-}(1.5,.5)(2.2,.5)
\psline[linewidth=1pt,linestyle=dotted]{-}(3.3,.5)(4,.5)
\psline[linewidth=1pt,linestyle=dotted]{-}(5.2,.5)(5.9,.5)
\psline[linewidth=1pt,linestyle=dotted]{-}(7,.5)(7.7,.5)
\psline[linewidth=1pt,linestyle=dotted]{-}(8.8,.5)(9.5,.5)
{\blue\uput[d](1.7,-.4){$_{E}$}\uput[d](3.7,-.4){$_{F}$}\uput[d](5.5,-.4){$_{G}$}\uput[d](7.3,-.4){$_{F^{-1}}$}
\uput[d](9.2,-.4){$_{H}$}}
{\red \uput[u](5.4,-.6){$_{\gz}$}}
\psline[linewidth=.5pt](4,1.5)(5.3,1.5)
{\small\uput[d](3,-1.4){$_{\t_{k_{t-1}}}$}
\uput[r](.8,0){$_{\t_{i_1}}$}\uput[l](2.7,0){$_{\t_{i_{t-1}}}$}\uput[r](3,0){$_{\t_{i_t}}$}
\uput[d](8.5,-1.4){$_{\t_{i_{t-1}}}$}
\uput[l](4.25,0){$_{\t_{i_s}}$}
\uput[r](4.45,0){$_{\t_{i_{s+1}}}$}
\uput[l](8,0){$_{\t_{i_t}}$}
\uput[r](8.3,0){$_{\t_{k_{t-1}}}$}
\uput[l](6.35,0){$_{\t_{i_{s+n}}}$}\uput[r](6.5,0){$_{\t_{i_s}}$}
\uput[r](9.8,0){$_{\t_{i_d}}$}}
\psline[linewidth=1pt,linestyle=dotted]{-}(1.1,-1.5)(1.9,-1.5)\psline[linewidth=1pt,linestyle=dotted]{-}(3.7,-1.5)(4.2,-1.5)
\psline[linewidth=1pt,linestyle=dotted]{-}(5.2,-1.5)(5.7,-1.5)\psline[linewidth=1pt,linestyle=dotted]{-}(7,-1.5)(7.5,-1.5)
\psline[linewidth=1pt,linestyle=dotted]{-}(9.2,-1.5)(9.6,-1.5)

\put(.85,-1.5){\circle*{.1}}\put(2.35,-1.5){\circle*{.1}}\put(3.35,-1.5){\circle*{.1}}\put(4.38,-1.5){\circle*{.1}}
\put(6.35,-1.5){\circle*{.1}}\put(7.85,-1.5){\circle*{.1}}\put(8.85,-1.5){\circle*{.1}}
\pscurve[linewidth=.5pt,linecolor=red]{-}(-.15,0)(2,-.5)(3,-.5)(7,-.5)(9,-.5)(10.85,0)
\psline[linewidth=.5pt,linecolor=blue]{->}(3.1,0)(2.6,0)\psline[linewidth=.5pt,linecolor=blue]{->}(4.1,0)(4.6,0)
\psline[linewidth=.5pt,linecolor=blue]{->}(6.1,0)(6.6,0)\psline[linewidth=.5pt,linecolor=blue]{->}(7.85,0)(8.35,0)
{\blue\uput[d](2.85,0.1){$_{\bz'}$}\uput[d](8.1,0.1){$_{\az'}$}
\uput[u](4.35,-0.1){$_{\az}$}\uput[u](6.35,-0.1){$_{\bz}$}}
\end{pspicture}
\end{center}
Hence $w(\ggz,\gz)$ can be given as follows:
{\small $$w(\ggz,\gz)=E-\t_{i_{t-1}}\overset{\bz'}{\leftarrow} \t_{i_t}-F-\t_{i_s}\overset{\az}{\ra} \t_{i_{s+1}}-G-\t_{i_{s+n}}\overset{\bz}{\ra} \t_{i_s}-F^{-1}-\t_{i_t}\overset{\az'}{\leftarrow} \t_{k_{t-1}}-H.$$}
We consider following two curves $\gz'$ and $\gz''$ in $(S,M)$:
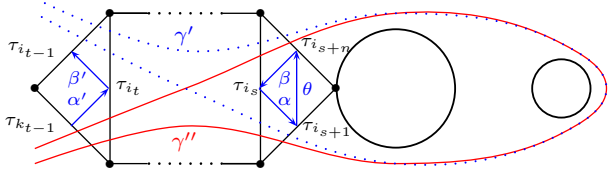
\begin{center}
\begin{pspicture}(0,-1.2)(9,1.2)
\pscircle[linewidth=.7pt](4.8,0){.8}    \pscircle[linewidth=.7pt](7,0){.4}
\put(4,0){\circle*{.1}}\put(3,1){\circle*{.1}}\put(1,1){\circle*{.1}}
\put(0,0){\circle*{.1}}\put(1,-1){\circle*{.1}}\put(3,-1){\circle*{.1}}
\psline[linewidth=.5pt]{-}(1.5,-1)(1,-1)(0,0)(1,1)(1.5,1)(1,1)(1,-1)\psline[linewidth=.5pt]{-}(2.5,1)(3,1)(4,0)(3,-1)(2.5,-1)(3,-1)(3,1)
\pscurve[linewidth=.5pt,linecolor=red]{-}(0,-.8)(2,0)(4.8,1)(7,.6)(7.6,0)(7,-0.6)(4.8,-1)(2,-0.5)(0,-1)
\psline[linewidth=1pt,linestyle=dotted]{-}(1.5,1)(2.5,1)\psline[linewidth=1pt,linestyle=dotted]{-}(1.5,-1)(2.5,-1)
\pscurve[linewidth=.8pt,linecolor=blue,linestyle=dotted]{-}(0,1.2)(2,0.5)(4.8,1)(7,.6)(7.6,0)(7,-0.6)(4.8,-1)(2,0)(0,1)
{\red\uput[d](2,-.4){$_{\gz''}$}}{\blue\uput[u](2,.4){$_{\gz'}$}}
{\small\uput[u](3.9,-.85){$_{\t_{i_{s+1}}}$}\uput[d](3.9,.85){$_{\t_{i_{s+n}}}$}
\uput[r](.9,0){$_{\t_{i_t}}$}
\uput[l](.5,.5){$_{\t_{i_{t-1}}}$}\uput[l](.5,-.5){$_{\t_{k_{t-1}}}$}
\uput[l](3.2,0){$_{\t_{i_s}}$}}

{\small\blue\uput[r](3.25,0){$_{\theta}$}
\uput[u](3.2,-.4){$_{\az}$}\uput[d](3.2,.4){$_{\bz}$}}
\psline[linewidth=.5pt,linecolor=blue]{->}(.35,-.5)(.85,0)
\psline[linewidth=.5pt,linecolor=blue]{->}(.85,0)(.35,0.5)\psline[linewidth=.5pt,linecolor=blue]{->}(3.35,.5)(2.85,0)
\psline[linewidth=.5pt,linecolor=blue]{<-}(3.35,-.5)(2.85,0)
\psline[linewidth=.5pt,linecolor=blue]{->}(3.35,-.5)(3.35,0.5)
{\small\blue\uput[u](.45,-.4){$_{\az'}$}\uput[d](.45,.4){$_{\bz'}$}}
\end{pspicture}
\end{center}
where $\gz'$ is given by the full curve and  $\gz''$ is described by the dotted curve. Then $w(\ggz,\gz')$ and $w(\ggz,\gz'')$ are given as follows
{\small
$$w(\ggz,\gz')=E-\t_{i_{t-1}}\overset{\bz'}{\leftarrow} \t_{i_t}-F-\t_{i_s}\overset{\az}{\ra} \t_{i_{s+1}}-G-\t_{i_{s+n}}\overset{\bz}{\ra}{\t_{i_s}}-F^{-1}-\t_{i_t}
\overset{\bz'}{\ra}\t_{i_{t-1}}-E^{-1}$$
$$w(\ggz,{\gz''})=H^{-1}-\t_{k_{t-1}}\overset{\az'}{\ra}\t_{i_t}-F-\t_{i_s}\overset{\az}{\ra} \t_{i_{s+1}}-G^{-1}-\t_{i_{s+n}}\overset{\bz}{\ra} \t_{i_s}-F^{-1}-\t_{i_t}\overset{\az'}{\leftarrow} \t_{k_{t-1}}-H.$$}

Hence there is a non-split exact sequence in $\mod J(Q_\ggz,W_\ggz)$£º
$$0\lra M(\ggz,\gz)\overset{{\tiny _{\left[
                              \begin{array}{c}
                                f \\
                                g \\
                              \end{array}
                            \right]
}}}{\lra}  M(\ggz,{\gz'})\oplus M(\ggz,{\gz''})\lra M(\ggz,\gz)\lra 0.$$
where $f$ sends the factor string $w_0=E-\t_{i_{t-1}}\overset{\bz'}{\leftarrow} \t_{i_t}-F-\t_{i_s}$ of $w(\ggz,\gz)$ to $w_0^{-1}$ as substring of $w(\ggz,\gz')$ and  $g$ identifies the factor string $w_1=\t_{i_t}-F-\t_{i_s}\overset{\az}{\ra} \t_{i_{s+1}}-G-\t_{i_{s+n}}\overset{\bz}{\ra} \t_{i_s}-F^{-1}-\t_{i_t}\overset{\az'}{\leftarrow} \t_{k_{t-1}}-H$ with $w_1^{-1}$ as substring of $w(\ggz,\gz'').$

\item[\bf (II.4)]$0<s<s+n=d.$

In this subcase, the subword $w(\ggz,\gz|_{[ r,r']})$ is a proper subword of $w(\ggz,\gz)$ which ends at the endpoint $e(\gz)$ of $\gz$.
This subcase is dual to subcase II.2.
\end{itemize}
\end{pf}
\medskip

\begin{rem}$\Ext^1_\mc(\gz,\gz)=0$ implies $\Ext^1_{J(Q_\ggz,W_\ggz)}(M(\ggz,\gz), M(\ggz,\gz))=0$ for any triangulation $\ggz$ of $(S,M)$. But $\Ext^1_{J(Q_\ggz,W_\ggz)}(M(\ggz,\gz), M(\ggz,\gz))=0$ for one triangulation $\ggz$ of $(S,M)$ does not imply $\Ext^1_\mc(\gz,\gz)=0$ in general. We reconsider $\da$ in Example \ref{ex2}. It is easy to see that
\begin{center}
\begin{pspicture}(1.5,-.5)(6.2,0.5)
\put(1.5,0){$_{M(\ggz,\da)=}$}
\put(3.1,0.3){$_{1}$}
\put(2.8,0){$_{~3~2~;}$}
\put(3.2,-.3){$_{~4}$}\put(3.3,-.6){$_{~5}$}

\put(4.5,0){$_{M(\ggz,\da[1])=}$}
\put(6.1,0.3){$_{1}$}
\put(6.15,0){$_{2~3.}$}
\put(6.2,-.3){$_{~4}$}
\end{pspicture}
\end{center}
Hence $\Ext^1_{J(Q_\ggz,W_\ggz)}(M(\ggz,\da),M(\ggz,\da))=0$ by the Auslander-Reiten formula in $\mod J(Q_\ggz,W_\ggz)$, since the non-zero morphism from $M(\ggz,\da)$ to $\t(M(\ggz,\da))$ $=M(\ggz,\da[1])$
factors through the injective module $I_5.$ But after applying a flip along 5, we get $\Ext^1_{\tilde{A}_4}(M(\ggz',\da),M(\ggz',\da))\neq 0$ which implies that $\Ext^1_\mc(\da,\da)\neq 0.$
\end{rem}
Similar as in Theorem \ref{intersection}, we can study the extensions of two different curves by their intersections in $(S,M).$

\begin{prop}Let $\gz$, $\da$ be two different curves in $(S,M)$, then $\Ext^1_\mc(\gz,\da)\neq 0 \neq\Ext^1_\mc(\da,\gz)$ if $I(\gz,\da)\neq 0$.
\end{prop}
\begin{pf}We only consider one of the intersections of $\gz$ and $\da$. After some flips (if needed), there is a triangulation $\ggz'$ such that $w(\ggz',\gz)$ and $w(\ggz',\da)$ share a common subword (associated to the intersection) $w=\t_{j_1}\cdots\t_{j_s}$ with $s\geq 1$, see the following picture
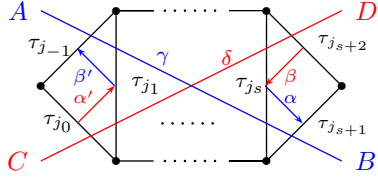
\begin{center}
\begin{pspicture}(-4,-1.1)(9,1)

\put(4,0){\circle*{.1}}\put(3,1){\circle*{.1}}\put(1,1){\circle*{.1}}
\put(0,0){\circle*{.1}}\put(1,-1){\circle*{.1}}\put(3,-1){\circle*{.1}}
\psline[linewidth=1pt,linestyle=dotted]{-}(1.6,-.5)(2.3,-.5)
\psline[linewidth=.5pt]{-}(1.5,-1)(1,-1)(0,0)(1,1)(1.5,1)(1,1)(1,-1)\psline[linewidth=.5pt]{-}(2.5,1)(3,1)(4,0)(3,-1)(2.5,-1)(3,-1)(3,1)
\pscurve[linewidth=.5pt,linecolor=red]{-}(0,-1)(2,0)(4,1)
\pscurve[linewidth=.5pt,linecolor=blue]{-}(4,-1)(2,0)(0,1)
{\small{\blue\uput[l](0,1){$A$}\uput[r](4,-1){$B$}}{\red\uput[l](0,-1){$C$}\uput[r](4,1){$D$}}}
{\red\uput[u](2.35,0.1){$_{\da}$}}{\blue\uput[u](1.5,0.1){$_{\gz}$}}
{\small\uput[u](3.9,-.85){$_{\t_{j_{s+1}}}$}\uput[d](3.9,.85){$_{\t_{j_{s+2}}}$}}
\uput[r](.9,0){$_{\t_{j_1}}$}
\psline[linewidth=1pt,linestyle=dotted]{-}(1.5,1)(2.2,1)\psline[linewidth=1pt,linestyle=dotted]{-}(1.5,-1)(2.2,-1)
\uput[l](.5,.5){$_{\t_{j_{-1}}}$}\uput[l](.5,-.5){$_{\t_{j_{0}}}$}
\uput[l](3.05,0){$_{\t_{j_s}}$}
\psline[linewidth=.5pt,linecolor=red]{->}(3.35,.5)(2.85,0)
\psline[linewidth=.5pt,linecolor=blue]{<-}(3.35,-.5)(2.85,0)
{\small\blue\uput[u](3.2,-.4){$_{\az}$}\uput[d](.45,.4){$_{\bz'}$}}
\psline[linewidth=.5pt,linecolor=red]{->}(.35,-.5)(.85,0)
\psline[linewidth=.5pt,linecolor=blue]{->}(.85,0)(.35,0.5)
{\small\red\uput[u](.45,-.4){$_{\az'}$}\uput[d](3.2,.4){$_{\bz}$}}
\end{pspicture}
\end{center}
where either $\t_{j_0}$ or $\t_{j_{-1}}$ might be a boundary arc, and one of $\t_{j_{s+1}}$ and $\t_{j_{s+2}}$ can also be a boundary arc. Hence
$w(\ggz',\gz)$ and $w(\ggz',\da)$ are of the form:
$$w(\ggz',\gz)=A-\t_{j_{-1}}\overset{\bz'}{\leftarrow}\t_{j_1}\leftrightsquigarrow\t_{j_s}\overset{\az}{\ra}\t_{j_{s+1}}-B,$$
$$w(\ggz',\da)=C-\t_{j_0}\overset{\az'}{\ra}\t_{j_1}\leftrightsquigarrow\t_{j_s}\overset{\bz}{\leftarrow}\t_{j_{s+2}}-D.$$
Consider the following two strings:
$$w(\ggz',\gz')=C-\t_{j_0}\overset{\az'}{\ra}\t_{j_1}\leftrightsquigarrow\t_{j_s}\overset{\az}{\ra}\t_{j_{s+1}}-B,$$
$$w(\ggz',\gz'')=A-\t_{j_{-1}}\overset{\bz'}{\leftarrow}\t_{j_1}\leftrightsquigarrow\t_{j_s}\overset{\bz}{\leftarrow}\t_{j_{s+2}}-D,$$
where $\gz'$ and $\gz''$ are given by following picture
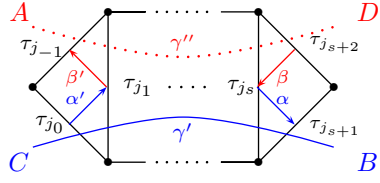
\begin{center}
\begin{pspicture}(-4,-1.1)(9,1)

\put(4,0){\circle*{.1}}\put(3,1){\circle*{.1}}\put(1,1){\circle*{.1}}
\put(0,0){\circle*{.1}}\put(1,-1){\circle*{.1}}\put(3,-1){\circle*{.1}}
\psline[linewidth=1pt,linestyle=dotted]{-}(1.8,0)(2.3,0)
\psline[linewidth=.5pt]{-}(1.5,-1)(1,-1)(0,0)(1,1)(1.5,1)(1,1)(1,-1)\psline[linewidth=.5pt]{-}(2.5,1)(3,1)(4,0)(3,-1)(2.5,-1)(3,-1)(3,1)
\pscurve[linewidth=1pt,linecolor=red,linestyle=dotted]{-}(4,.8)(2,.4)(0,.8)
\pscurve[linewidth=.5pt,linecolor=blue]{-}(4,-.8)(2,-.4)(0,-.8)
{\small\red\uput[u](2,0.3){$_{\gz''}$}}{\blue\uput[d](2,-0.3){$_{\gz'}$}}
{\small{\blue\uput[l](0,-1){$C$}\uput[r](4,-1){$B$}}{\red\uput[l](0,1){$A$}\uput[r](4,1){$D$}}}
{\small\uput[u](3.9,-.85){$_{\t_{j_{s+1}}}$}\uput[d](3.9,.85){$_{\t_{j_{s+2}}}$}}
\uput[r](.9,0){$_{\t_{j_1}}$}
\psline[linewidth=1pt,linestyle=dotted]{-}(1.5,1)(2.2,1)\psline[linewidth=1pt,linestyle=dotted]{-}(1.5,-1)(2.2,-1)
\uput[l](.5,.5){$_{\t_{j_{-1}}}$}\uput[l](.5,-.5){$_{\t_{j_{0}}}$}
\uput[l](3.05,0){$_{\t_{j_s}}$}
\psline[linewidth=.5pt,linecolor=red]{->}(3.35,.5)(2.85,0)
\psline[linewidth=.5pt,linecolor=blue]{<-}(3.35,-.5)(2.85,0)
{\small\blue\uput[u](3.2,-.4){$_{\az}$}\uput[u](.45,-.4){$_{\az'}$}}
\psline[linewidth=.5pt,linecolor=blue]{->}(.35,-.5)(.85,0)
\psline[linewidth=.5pt,linecolor=red]{->}(.85,0)(.35,0.5)
{\small\red\uput[d](3.2,.4){$_{\bz}$}\uput[d](.45,.4){$_{\bz'}$}}
\end{pspicture}
\end{center}
then there is a non-split exact sequence in $\mod J(Q_{\ggz'},W_{\ggz'})$£º
$$0\lra M(\ggz',\gz)\overset{{\tiny _{\left[
                              \begin{array}{c}
                                f \\
                                g \\
                              \end{array}
                            \right]
}}}{\lra}  M(\ggz',{\gz'})\oplus M(\ggz',{\gz''})\lra M(\ggz',\da)\lra 0.$$
where $f$ sends the factor string $w_0=\t_{j_1}\leftrightsquigarrow\t_{j_s}\overset{\az}{\ra}\t_{j_{s+1}}-B$ of $w(\ggz',\gz)$ to $w_0$ as substring of $w(\ggz',\gz')$ and  $g$ identifies the factor string $w_1=A-\t_{j_{-1}}\overset{\bz'}{\leftarrow}\t_{j_1}\leftrightsquigarrow\t_{j_s}$ as substring of $w(\ggz',\gz'').$ Therefore
$$\Ext^1_{J(Q_{\ggz'},W_{\ggz'})}(M(\ggz',\da),M(\ggz',\gz))\neq 0$$ which implies that
$\Ext^1_\mc(\da,\gz)\neq 0,$ and since $\mc$ is 2-Calabi-Yau, we have $\Ext^1_\mc(\gz,\da)\cong D\Ext^1_\mc(\da,\gz)\neq 0$
\end{pf}
The above proposition and Corollary \ref{cor1} imply:
\begin{cor}Let $\gz$, $\da$ be two different internal arcs in $(S,M)$, then
$\Ext^1_\mc(\gz,\da)=0$ if and only if $I(\gz,\da)=0.$
\end{cor}

\begin{cor}There is a bijection between triangulations of $(S,M)$ and cluster-tilting objects of $\mc_{(S,M)}$. In particular, each indecomposable object without self-extensions is reachable from the cluster-tilting object $T_\ggz$ (the initial cluster-tilting object).
\end{cor}
\begin{pf}Corollary \ref{cor1} implies that each triangulation of $(S,M)$ yields
a cluster-tilting object
in $C_{(S,M)}$. Hence it suffices to prove that each cluster tilting object
$T=T_1\oplus\cdots\oplus T_n$ gives a triangulation of $(S,M)$. By Theorem
\ref{intersection}, we can assume $T=\t_1\oplus\cdots\oplus\t_n$ where each
$\t_i$ is an internal arc in $(S,M)$ corresponding to $T_i.$ The definition of
a triangulation and the above corollary yield a unique triangulation
$\ggz_T=\{\t_1,\cdots,\t_n,\t_{n+1},\cdots,\t_{n+m}\}$ where
$\t_{n+1},\cdots,\t_{n+m}$ are boundary arcs in $(S,M).$
\end{pf}

\end{document}